\documentclass[a4paper,12pt,leqno]{amsart}
\usepackage{hyperref}
\usepackage{amsmath}
\usepackage{amssymb}
\usepackage{amsthm}
\usepackage{amsopn}
\usepackage{amscd}
\usepackage{fullpage}
\usepackage[utf8]{inputenc}
\usepackage{graphics, xcolor}
\usepackage{textcomp} 
\usepackage[all]{xy}
\usepackage{lscape}

\usepackage{verbatim}
\usepackage{mathrsfs}

\newcommand{\scrC}{\mathscr{C}}
\newcommand{\scrD}{\mathscr{D}}
\newcommand{\scrE}{\mathscr{E}}
\newcommand{\scrF}{\mathscr{F}}

\newcommand{\scrH}{\mathscr{H}}
\newcommand{\scrI}{\mathscr{I}}

\newcommand{\scrK}{\mathscr{K}}

\newcommand{\scrM}{\mathscr{M}}

\newcommand{\scrR}{\mathscr{R}}

\newcommand{\frZ}{\mathfrak{Z}}

\newcommand{\fra}{\mathfrak{a}}
\newcommand{\frb}{\mathfrak{b}}

\newcommand{\frg}{\mathfrak{g}}
\newcommand{\frh}{\mathfrak{h}}

\newcommand{\frk}{\mathfrak{k}}

\newcommand{\fro}{\mathfrak{o}}
\newcommand{\frp}{\mathfrak{p}}

\newcommand{\frs}{\mathfrak{s}}
\newcommand{\frt}{\mathfrak{t}}
\newcommand{\fru}{\mathfrak{u}}

\newcommand{\bbC}{\mathbb{C}}

\newcommand{\bbN}{\mathbb{N}}

\newcommand{\bbR}{\mathbb{R}}

\newcommand{\bbZ}{\mathbb{Z}}

\newcommand{\caA}{\mathcal{A}}

\newcommand{\caC}{\mathcal{C}}
\newcommand{\caD}{\mathcal{D}}

\newcommand{\caF}{\mathcal{F}}

\newcommand{\caL}{\mathcal{L}}
\newcommand{\caM}{\mathcal{M}}

\newcommand{\sfH}{\mathsf{H}}
\newcommand{\sfI}{\mathsf{I}}

\newcommand{\coker}{\mathrm{coker} }

\theoremstyle{plain}
\newtheorem{thm}{Theorem}[section]
\newtheorem{lemma}[thm]{Lemma}
\newtheorem{cor}[thm]{Corollary}
\newtheorem{prop}[thm]{Proposition}
\newtheorem{conj}[thm]{Connjecture}
\theoremstyle{definition}
 \newtheorem{defi}[thm]{Definition}

\theoremstyle{definition}

\newtheorem{rmk}[thm]{Remark}
\newtheorem{example}[thm]{Example}

\newcommand{\Ad}{\mathrm{Ad}}
\newcommand{\ad}{\mathrm{ad}}

\newcommand{\End}{\mathrm{End}}
\newcommand{\Hom}{\mathrm{Hom}}
\newcommand{\im}{\mathrm{Im}\, }

\newcommand{\Id}{\mathrm{Id}}

\newcommand{\ind}{\mathrm{ind}}

\newcommand{\SO}{\mathbf{SO}}

\newcommand{\Spin}{\mathbf{Spin}}
\newcommand{\Cl}{\mathrm{Cl}}
\newcommand{\Ext}{\mathrm{Ext}}

\newcommand{\Rk}{\mathrm{Rk}}

\newcommand{\Ep}{\mathrm{EP}}

\newcommand{\bil}[2]{\langle  #1,#2 \rangle }
\newcommand{\bilo}{\langle  .\, , . \rangle }

\newcommand{\Hdir}{\sfH_{\mathrm{Dir}}}

\newcommand{\Idir}{\sfI_{\mathrm{Dir}}}

\def \proof {\noindent \underline{\sl Proof}. }

\newcommand{\EP}{\mathbf{EP}}

\begin{document}

\numberwithin{equation}{section}

\title{Euler-Poincaré pairing, Dirac index and
 elliptic pairing for Harish-Chandra modules}

\author{ D. RENARD}
\address{Centre de math\'ematiques Laurent Schwartz, Ecole Polytechnique, 91128 Palaiseau Cedex}
\email{renard@math.polytechnique.fr}
\date{\today}

\begin{abstract}Let $G$ be a connected real reductive group with  maximal 
compact subgroup $K$ of equal rank, and let $\scrM$ be the category of Harish-Chandra modules for $G$. 
We relate  three differentely defined pairings between two finite length  modules $X$ and $Y$
 in $\scrM$   :  the Euler-Poincaré pairing,  the natural pairing between the Dirac indices  of $X$ and $Y$, 
 and the elliptic pairing of \cite{Ar}.  
 (The Dirac index $\Idir(X)$ is a virtual finite dimensional representation of $\widetilde K$,   the spin  double cover of $K$.)

Analogy with the case of Hecke algebras studied in \cite{CT} and \cite{COT} and a formal (but not 
rigorous) computation lead us to conjecture that the first two pairings coincide.
 In the second part of the paper, we show that they are both computed as the indices of   Fredholm pairs
 (defined here in an algebraic sense) of operators acting on the same spaces.
 
We construct index functions $f_X$ for any  finite length  Harish-Chandra 
module $X$. These functions are very cuspidal in the sense of Labesse, and their 
orbital integrals on elliptic elements coincide with the character of $X$. From this we deduce that
 the Dirac index  pairing coincide with the elliptic pairing.

  These results are  the archimedean   analog of results of Schneider-Stuhler \cite{ScSt} for $p$-adic groups.
\end{abstract}

\maketitle

\section{Introduction}

Our goal is to establish the archimedean   analog of results of Schneider-Stuhler \cite{ScSt}.
Let us first describe these results.
Let $G$ be the group of rational points of a connected algebraic reductive group defined
over a local non archimedean field $F$ of characteristic $0$.
 Assume that $G$ has compact center.  Let $\scrM=\scrM(G)$ be the category of smooth
representations of $G$.  This is also the category of non-degenerate modules over the
Hecke algebra $\scrH=\scrH(G)$ of $G$. It is known from the work of J. Bernstein that $\scrM$
has finite cohomological dimension. Furthermore, for any finitely generated  module $(\pi,V)$ in $\scrM$, 
 Schneider and Stuhler    construct an explicit  resolution
 of $(\pi,V)$  by finitely generated projective modules. 
 They also establish a general theory of Euler-Poincaré
 functions for modules of finite length, generalizing results of Kottwitz (\cite{Kot}).
 Namely, for any finite length modules $(\pi,V)$, $(\pi',V')$ in $\scrM$, 
 one can define their Euler-Poincaré pairing :
 \begin{equation}\label{EPp}  \EP(\pi,\pi')=\sum_i (-1)^i \; \dim \Ext_{\scrM}^i(\pi, \pi').
   \end{equation} 
 They construct functions $f_\pi$ (Euler-Poincaré functions) in $\scrH$ such that 
 \begin{equation}\label{EPfunc}    \EP(\pi,\pi')= \Theta_{\pi'}(f_\pi) \end{equation}
 where $\Theta_{\pi'}$ is the distribution-character of $\pi'$ (a linear form on $\scrH$).
 Following \cite{Dat} and \cite{Vig}, let us now give another point of view on  these functions.
 
 Let $\scrK (G)$ be the Grothendieck group of finitely generated projective modules in $\scrM$. 
 Since $\scrM $ has finite cohomological dimension, this is also the Grothendieck group of all 
 finitely generated modules in $\scrM$. 
  Let $\scrR=\scrR(G)$ be the Grothendieck
  group of finite length modules in $\scrM$. If $(\pi,V)$ is a finitely generated  (resp. finite length) 
  module  in $\scrM$, we denote by $[\pi]$ its image in  Grothendieck group  $\scrK$ (resp. $\scrR$).
   Set    $\scrK_\bbC=\scrK\otimes_\bbZ\bbC$ and $\scrR_\bbC=\scrR\otimes_\bbZ\bbC$. 
   Let  $(\pi,V)$ be a finite length module in $\scrM$, and 
  \[\cdots  \longrightarrow P_{i+1} \longrightarrow P_{i}  \longrightarrow \cdots \longrightarrow P_{1}
 \longrightarrow V  \longrightarrow 0\]
be a resolution of $\pi$ by finitely generated projective modules. Then 
 \begin{equation}\label{EPmap}  \Ep : \; \scrR \longrightarrow  \scrK, \quad 
 [\pi] \mapsto \sum_i (-1)^i \; [P_i]  
   \end{equation} 
is a well-defined map.  
 
Let $\overline {\scrH}=\overline {\scrH(G)}=\scrH/[\scrH,\scrH] $ be the abelianization   
 of $\scrH$. The Hattori rank map 
  \begin{equation}\label{HatorriRk}  \Rk : \; \scrK \longrightarrow  , \quad \overline {\scrH} 
   \end{equation} 
 is defined as follows. Let $P$ be a  finitely generated projective module in $\scrM$. Write 
 $P$ as a direct factor of some $\scrH^n$ and let $e\in \End_\scrH(\scrH^n)$ be the projector
 onto $P$. Then the trace of $e$ is an element of $\scrH$, and its image in $\overline {\scrH}$
is well-defined.  This defines $\Rk([P])$. An alternative description of $\overline {\scrH}$ 
as the "cocenter" of the category $\scrM$ gives a natural definition of $\Rk$ for 
 finitely generated projective modules  (see \cite{Dat}, \S 1.3).

 Let $\caD'(G)$ be the space of distributions on $G$, and 
 let $\caD'(G)^G$ be the subspace of invariant distributions. 
 Fix a Haar measure on $G$, so  that $\scrH$ is identified with 
  the convolution algebra of compactly supported smooth functions
on $G$. The orthogonal of $\caD'(G)^G$ in $\scrH$
for the  natural pairing between $\caD'(G)$ and $\scrH$ is exactly $[\scrH,\scrH]$, so there is 
an induced non-degenerate  pairing :
\[  \caD'(G)^G \times \overline{\scrH}  \rightarrow \bbC, \quad (T,f) \mapsto \bil{T}{f}=T(f). \]
Let $(\pi',V')$ be a finite length module in $\scrM$. Its distribution character 
$\Theta_{\pi'}$ is an element of   $\caD'(G)^G$. This defines a pairing 
  \begin{equation}\label{PairingRH}   \scrR_\bbC \times \overline{\scrH} \rightarrow \bbC, \quad
  (\pi',f)\mapsto \bil{\Theta_{\pi'}}{f} = \Theta_{\pi'}(f) 
   \end{equation} 
 
With the notation above,  we have the following identity : for all
finite length modules $(\pi,V),(\pi',V')$ in $\scrM$, 
 \begin{equation}\label{EPpform}  \EP(\pi,\pi')=\bil{\Theta_{\pi'}}{\Rk\circ \Ep(\pi)}.
   \end{equation} 
 Thus, the image in $\overline {\scrH}$ of the Euler-Poincaré function 
 $f_\pi$ constructed by Schneider and Stuhler  is $\Rk\circ \Ep(\pi)$ (\cite{Dat}, lemma 3.7).

 There is a third way of seeing  the space $\overline {\scrH}$, namely, as the space of orbital
  integrals on $G$. More precisely, recall that for a regular semisimple element $x$  in $G$, one can 
  define the orbital integral of $f \in \scrH$ at $x$ as 
  \[ \Phi(f,x)= \int_{G/T}  f(gxg^{-1}) \; d\dot{g} ,  \]
  where $T$ is the unique maximal torus containing $x$, and  $d\dot{g}$  
 is an invariant measure on $G/T$.
When $f$ is fixed, $x \mapsto  \Phi(f,x)$ is a smooth invariant function on $G_{\mathrm{reg}}$, and we denote  
by  $\scrI(G)$ the image  of $\Phi :\, f \in \scrH \mapsto \Phi(f,.) $ in the space of  smooth invariant functions on 
$G_{\mathrm{reg}}$. This space  can be explicitely described (by properties of orbital integrals, see \cite{Vig2}). 
Furthermore the kernel of $\Phi$ is exactly $[\scrH,\scrH]$ (this is called the geometric density theorem, 
{\sl i.e.} the density of the space generated by the distributions $f \mapsto \Phi(f,x)$, $ x \in 
G_{\mathrm{reg}} $, in $\caD'(G)^G$).
Thus, we have an exact sequence 
\[ 0\longrightarrow [\scrH,\scrH]\longrightarrow \scrH  \longrightarrow  \scrI(G)\longrightarrow
0\]
and $\overline{\scrH}$ is identified with the space  $\scrI(G)$ of orbital integrals.

Let us denote by $G_{\mathrm{ell}}$ the space of regular semisimple  elliptic elements in 
$G$ ({\sl i.e.} elements whose centralizer is a compact maximal  torus in $G$), and by 
$\scrH_c$ the subspace of functions $f \in \scrH$ such that $\Phi(f,x)=0$
when $x$ is not elliptic. Acccordingly, we write $\overline{\scrH}_c$
for the image of this subspace in $\overline{\scrH}$. This latter space is isomorphic to the 
subspace   $\scrI(G)_c$ of orbital integrals which vanish outside  $G_{\mathrm{ell}}$.
The Selberg principle for Euler-Poincaré functions (\cite{ScSt}, Rmk II.4.11) asserts that 
the functions $f_\pi$ are in $\scrH_c$. In fact, it is a theorem of Brylinski-Blanc
\cite{BlBr}, (see also \cite{Dat} and \cite{Dat2}) that the image Hattori rank map $\Rk$  is exactly 
$\overline{\scrH}_c$.

 Recall that by Harish-Chandra regularity theorem, the character $\Theta_{\pi'}$  of a finite length 
 module  $(\pi',V')\in \scrM$ is given by a locally integrable function,  denoted by $\theta_{\pi'}$, 
 on $G$, which is smooth on the open dense subset of regular semisimple elements $G_{\mathrm{reg}}$, 
 {\sl i.e}  for all $f \in \scrH$, 
 \[  \bil{\Theta_{\pi'}}{f}=\int_{G_{\mathrm{reg}}} \theta_{\pi'}(x)\, f(x) \; dx.   \]
 (Once again, we have chosen a Haar measure on $G$).

Following Kazhdan \cite{Kaz}, let us denote by $C^{\mathrm{ell}}$ the set of regular semisimple elliptic
conjugacy classes in $G$. Then any orbital integral $\Phi(f,.)$ or any character function $\theta_{\pi'}$
can be viewed as  functions on $C^{\mathrm{ell}}$.
By \cite{Kaz}, \S3, Lemma 1, there is a unique measure $dc$ on $C^{\mathrm{ell}}$
such that for all $f \in \scrH$ with support in $G_{\mathrm{ell}}$, 
\[ \int_G  f(x) \; dx= \int_{C^{\mathrm{ell}}} \Phi(f,c)\; dc. \]
The elliptic pairing between two  finite length 
 modules  $(\pi,V)$ and $(\pi',V')$ in  $\scrM$ is then defined by 
 \begin{equation} \label{EllPairing} 
 \bil{\pi}{\pi'}_{\mathrm{ell}}= \int_{C^{\mathrm{ell}}} \theta_\pi(c) \, \theta_{\pi'}(c^{-1})\; dc.
  \end{equation}
  Let us now relate this elliptic pairing to the Euler-Poincaré pairing (\ref{EPp}).
On  a regular semisimple  elliptic element $x$, the orbital integrals of the Euler-Poincaré functions 
$f_\pi$ at $x$ coincide with the character  $\theta_{\pi}$ of $\pi$ at $x^{-1}$ (\cite{ScSt} Thm III.4.16) : 
\begin{equation}\label{psCoeff} (\forall x \in   G_{\mathrm{ell}})   , \qquad   \Phi(f_\pi,x)= \theta_{\pi}(x^{-1}).
\end{equation} 
 Therefore, the following formula,  which was conjectured by Kazhdan \cite{Kaz} is established 
 (\cite{ScSt} Thm III.4.21): 
\begin{equation}\label{EPelliptic}  \EP(\pi,\pi')=\bil{\Theta_{\pi'} } {f_\pi }=\int_{G_{\mathrm{ell}}}
\theta_{\pi'}(x) \, f_\pi(x) \; dx =\int_{C^{\mathrm{ell}}}
\theta_{\pi'}(c) \, \theta_{\pi}(c^{-1}) \; dc=  \bil{\pi}{\pi'}_{\mathrm{ell}}.
   \end{equation} 
Let us mention that it  is explained in \cite{Dat} how these results can be formulated and 
proved when the center of  $G$ is not compact.

\medskip

Many of  the objects and  results described above make sense over an Archimedean field as well. 
So assume now that $G$ be the group of rational points of a connected algebraic reductive group defined
over $\bbR$ and, for simplicity,  that $G$ is connected (as a Lie group), with compact center.
A little bit more generally, we  may assume without introducing technical complications that $G$
 is actually a finite covering 
of such a group. Fix a maximal compact subgroup  $K$ of $G$.
In section \ref{psCoeff}, also for simplicity, we will assume that $G$ and $K$
have same rank, {\sl i.e.} $G$ admits discrete series representations (this implies the 
compact center condition)\footnote{In a future version of this paper, we hope to get rid of some of these
non necessary assumptions}. The category of representations 
we are now considering is  the category of Harish-Chandra modules $\scrM=\scrM(\frg,K)$ (\cite{KV}).
It is also known that $\scrM$ has finite cohomological dimension, and furthermore
 resolutions of finite length modules by explicit finitely generated projective
modules can be constructed (\cite{KV}, \S II.7). Thus, the Euler-Poincaré
pairing of two modules of finite length  is well-defined by (\ref{EPp}).
The elliptic pairing (\ref{EPelliptic}) is also defined, and the purpose of this paper is to relate
these two pairings by constructing  Euler-Poincaré functions  in that case, 
{\sl i.e.} compactly supported smooth functions on $G$ satisfying (\ref{EPfunc}),
 but also  (\ref{psCoeff}), so that the formula  (\ref{EPelliptic}) is valid. 
   An interesting aspect  is that this  relation involves a third pairing defined through  Dirac cohomology.  
 
 Dirac cohomology of Harish-Chandra modules was  introduced by D. Vogan (see \cite{HPbook}
  for an exposition of this theory). Dirac cohomology of a Harish-Chandra module $X$ of finite length 
  consists in  a finite   dimensional representation $\Hdir(X)$  of the compact group $\widetilde K$, the spin cover of $K$.
In fact, $\Hdir(X)$ is $\bbZ /2\bbZ$-graded, and a slightly more refined invariant, the Dirac index 
$\Idir(X)$ can be defined as the formal difference $\Idir(X)=\Hdir(X)^{\bar 0}-\Hdir(X)^{\bar 1}$
between the even and odd part of $\Hdir(X)$,
 a virtual   finite  dimensional representation of $\widetilde K$.
If $\left[\sigma,\sigma'\right]_{\widetilde K}$ is the usual pairing between two  virtual 
finite  dimensional representations $\sigma, \sigma'$  of $\widetilde K$, then 
the third pairing between two Harish-Chandra modules $X$ and $Y$  of finite length that we introduce  is 
\begin{equation}\label{DirPair}  \bil{X}{Y}_{\mathrm{Dir}}= \left[ \Idir(X),\Idir(Y)\right]_{\widetilde K}.
\end{equation}
A formal (but not rigorous)  computation given in Section \ref{EPdir} of the paper
lead us to conjecture that 
\begin{equation}\label{ConjEPD}  \EP(X,Y)=\bil{X}{Y}_{\mathrm{Dir}}, \end{equation}
for any Harish-Chandra modules $X$, $Y$ of finite length. 
Interestingly enough, to give a complete  proof is harder than it seems first, and we did not  succeed
so far. In Section \ref{Proof}  we  show that both sides of (\ref{ConjEPD}) are  indices of suitable  Fredholm
pairs defined on the same spaces. We found the definition of Fredholm pairs  (and some properties of their indices) in 
   \cite{Amb} and we adapted it to our purely algebraic setting. Let us just say here that the index of 
  a Fredholm pair is a generalization of the index of a Fredholm operator. The relevant material is exposed in an appendix.
The first of these Fredholm pair is given by the complex computing the $\Ext$ groups and its differential, and its index is 
the Euler-Poincaré characteristic of the complex, {i.e.} the left-hand side of (\ref{ConjEPD}). 
The second Fredholm pair is given (on the same space) by actions of Dirac operators and its index is the 
right-hand side of (\ref{ConjEPD}). The conjecture thus boils down to the equality of these two indices. Of course this kind
of result looks familiar. This gives some credit to the conjecture and points out a direction to look for a proof.

In section \ref{EPLab}, we use Labesse  "index" functions (\cite{Lab})
 to construct index functions \footnote{Assuming the conjecture, we could call these Euler-Poincaré functions.}
  $f_X$ for any finite length Harish-Chandra 
 module $X$, {\sl i.e.} $f_X$ is a smooth compactly supported function on $G$ satisfying 
\[ \bil{X}{Y}_{\mathrm{Dir}},= \Theta_Y(f_X) \]
for any finite length Harish-Chandra module $Y$. 
These functions are very cuspidal in the sense
of Labesse, thus their orbital integrals $\Phi(f_X,x)$ vanish on non elliptic elements.  
In section \ref{OrbInt}, we show that when $x$ is an elliptic element in $G$,
the formula (\ref{psCoeff}) is valid in this context. 
 This easily  implies that  the elliptic pairing coincide with the Dirac pairing.
  The proof uses  the density of tempered characters in the space of invariant distributions to reduce
the problem to the case of limits of discrete series. In the  case of discrete series, the relevant results are well-known
and their extension to the case of limits of discrete series is relatively easy.


\bigskip 

The idea that the elliptic pairing for Harish-Chandra modules is related to Dirac index
originates in the papers \cite{CT} and \cite{COT} where the result is  established for Hecke algebra.
The role of the group $\widetilde K$ is played in that context by the spin cover 
$\widetilde W$ of the relevant Weyl group. Since modules for the Hecke algebra are finite dimensional, the difficulties we met in 
proving (\ref{ConjEPD}) do not appear.

\bigskip 

We learned from G. Zuckerman that   
he  obtained  results on Euler-Poincaré pairing for Harish-Chandra modules 
in the late 70's, but that these were never published. Pavle Pand\v zi\'c and Peter Trapa   informed us that they 
were also led to  the identity (\ref{ConjEPD})
in a work in progress with Dan Barbasch.

\section{Dirac cohomology and Dirac index of $(\frg,K)$-modules}\label{Dgk}

\subsection{Notation and structural facts}\label{gK}
Let  $G$ be  a connected real reductive Lie group with Cartan involution $\theta$ such that $K=G^\theta$ 
is a maximal compact subgroup of $G$. Let us denote by 
$\frg_0$  the Lie algebra of $G$, $\frg$ its complexification,   with  Cartan involutions 
also denoted by $\theta$.

We fix an  invariant nondegenerate symmetric bilinear form $B$ on $\frg_0$, 
 extending the Killing form  on the semisimple part of $\frg_0$. Let 
\[ \frg_0=\frk_0 \stackrel{\theta}{\oplus}\frp_0, \qquad \frg=\frk \stackrel{\theta}{\oplus}\frp\]
be the  Cartan decompositions of $\frg_0$ and $\frg$ respectively.
We assume furthermore that in extending the Killing form, we  made sure that $B_{|\frp_0}$
remains definite positive and $B_{|\frk_0}$ definite negative.
We assume for simplicity that the center of $G$ is compact. This implies that the dimension of $\frp_0$
is even.

Let $\Cl(\frp)=\Cl(\frp;B)$ be  the Clifford algebra of $\frp$ with respect to $B$.
 The Clifford algebra  $\Cl(\frp)$  is the quotient of the tensor algebra $T(\frp)$ by the
  two-sided ideal generated by elements of the form
\[v\otimes w+w\otimes v+2B(v,w)\, 1, \qquad (v,w \in V).\]
Some authors use relations differing by a sign.   We follow here the convention of  \cite{HPbook}.
We refer to \cite{HPbook}, \cite{Kos1} or \cite{Mein} for material on Clifford algebras and spinor
modules.

 Notice that all the terms in the above expression are  of even degrees ($2$ or $0$). 
 Thus,  the  graded algebra structure on $T(\frp)$ induces a  filtered algebra structure on 
 $\Cl(\frp)$, but also a structure of $\bbZ/2\bbZ$-graded algebra  ({\sl i.e.}  
 a super algebra structure).
 Simply put, the $\bbZ/2\bbZ$-grading and the filtration
  are defined by the condition that the generators
 $X\in \frp$ of  $\Cl(\frp)$ are odd, of filtration  degree $1$.
   We get a decomposition
 \[  \Cl(\frp)= \Cl^{\bar 0}(\frp) \oplus \Cl^{\bar 1}(\frp). \]
 Recall also that the Clifford algebra $\Cl(\frp)$ is isomorphic as a  $\bbZ/2\bbZ$-graded 
 vector space to  the exterior algebra $\bigwedge \frp$ by the Chevalley isomorphism. It is convenient 
 to identify  $\Cl(\frp)$ to  $\bigwedge \frp$, and to see the latter as a  $\bbZ/2\bbZ$-graded 
 vector space   endowed with two different algebra products, the graded commutative  wedge product 
 (denoted $x\wedge y$)  and the Clifford product (denoted $xy$).
 Let us denote be $\SO(\frp_0)$ (resp. $\SO(\frp)$) the special orthogonal group of $(\frp_0,B)$ 
 (resp. $(\frp, B)$).  The subspace $\bigwedge^2\frp$ is stable under the Clifford Lie bracket
$[x,y]_\Cl=xy-yx$ on $\bigwedge \frp=\Cl(\frp)$ and is isomorphic as a Lie algebra to the Lie algebra
$\frs\fro(\frp)$ of   $\SO(\frp)$.

 We denote by $\widetilde K$ the spin double cover of $K$, i.e., the pull-back of the covering 
map $\mathbf{Spin}(\frp_0)  \rightarrow \mathbf{SO}(\frp_0)$   by  the adjoint action map 
$\Ad_{|\frp_0} :\, K \rightarrow \SO(\frp_0)$.
The compact groups $\Spin(\frp_0)$ and   $\SO(\frp_0)$ embed in their complexification 
 $\Spin(\frp)$ and   $\SO(\frp)$, so we get the following diagram : 
\[ \xymatrix{  
  \widetilde{K} \ar[r]  \ar[d]    &    \Spin(\frp_0)  \ar@{^{(}->}[r]  \ar[d]  &   
   \Spin(\frp)\ar@{^{(}->}[r]  \ar[d]  
 & \Cl^{\bar 0}(\frp)^\times  \\
   K \ar[r]^{\Ad_{|\frp_0}}  & \SO(\frp_0) \ar@{^{(}->}[r]  & \SO(\frp)  & \\
}\]
The complexification of the  differential at the identity of the Lie 
group morphism  $\Ad_{|\frp_0} :\, K \rightarrow \SO(\frp_0)$, is the Lie algebra
morphism
\begin{equation*}
 \ad_{|\frp}:\, \frk \rightarrow \frs\fro(\frp), \quad X \mapsto \ad(X)_{|\frp}   \end{equation*}
Let us denote by $\alpha$ the composition of this map with the identification 
between $\frs\fro(\frp)$ and $\bigwedge^ 2 \frp$ and the inclusion of  $\bigwedge^ 2 \frp$ in  $\bigwedge \frp=\Cl(\frp)$ : 
\begin{equation}\label{alpha1}
\alpha :\, \frk \rightarrow \Cl(\frp).
\end{equation}

A key role is  played in the theory of Dirac cohomology of Harish-Chandra modules
 by the associative  $\bbZ/2\bbZ$-graded superalgebra 
$\caA=U(\frg)\otimes \Cl(\frp)$. The    $\bbZ/2\bbZ$-grading comes from the  $\bbZ/2\bbZ$-grading on 
$\Cl(\frp)$, {\sl i.e.} elements in  $U(\frg)\otimes 1$  are even. The 
 super Lie algebra structure  on  $\caA$ is given by the (super)commutator bracket $[.\, ,.]_\caA$.

The group $K$ acts on $U(\frg)$ through $K\subset G$ by the adjoint action, and 
on $\Cl(\frp)$ through the map  $\widetilde K \rightarrow \Cl^{\bar 0}(\frp)^\times $ 
in  the first row of the  diagram above and conjugation in $\Cl(\frp)$ (this action of $\widetilde K$
 on $\Cl(\frp)$ factors through $K$).
 Thus we get a linear action of  $K$ on $\caA$.
Differentiating  this action at the identity, and taking the complexification,  we get
a Lie algebra representation of $\frk$ in $\caA$.
This represention can be described as follows.  The map  (\ref{alpha1}) is used to define
a map 
\begin{equation*}
\Delta:\, \frk \longrightarrow  \caA= U(\frg)\otimes \Cl(\frp), \quad 
  \Delta(X)=X\otimes 1+1\otimes\alpha(X)
\end{equation*}
which is a morphism of Lie algebra (it takes values in the even part 
of the super Lie algebra $\caA$). Thus it extends to an algebra morphism
\begin{equation}\label{Delta}
\Delta:\, U(\frk) \longrightarrow  \caA= U(\frg)\otimes \Cl(\frp).
\end{equation}
The action of an element $X \in \frk$ on $\caA$ is then given by 
the adjoint action of $\Delta(X)$, {\sl i.e.} $a\in \caA \mapsto [\Delta(X),a]_\caA$.
We denote by $\caA^K$ (resp. $\caA^\frk$) the subalgebra of $K$-invariants (resp. 
$\frk$-invariants)  in $\caA$. Since $K$ is assumed to be connected, $\caA^K=\caA^\frk$.

Let us now recall some facts about Clifford modules. 
\begin{thm}

 Suppose that $n=\dim(\frp)$ is even. Then there are : 
\smallskip 

$\bullet$ two isomorphism classes of irreducible $\bbZ/2\bbZ$-graded $\Cl(\frp)$-modules, 

$\bullet$ one  isomorphism class of irreducible ungraded $\Cl(\frp)$-modules, 

$\bullet$ two isomorphism classes of irreducible  $\Cl^{\bar 0}(\frp)$-modules. 
\medskip





\medskip

\end{thm}

In the case  $\dim(\frp)$ is even, an irreducible  ungraded $\Cl(\frp)$-module $S$ can be realized as follows :  
choose a decomposition $\frp=U\oplus U^*$ into dual isotropic subspaces. 
As the notation indicates, $U^*$ is identified to the dual of $U$ by the bilinear form $B$.
Set  $S=\bigwedge U$. Let $U$ act on $S$ by wedging and $U^*$ by contracting. 
The decomposition $\textstyle \bigwedge U = \bigwedge^{\bar{0}} U\oplus  
 \bigwedge^{\bar{1}}U$ give a decomposition $S=S^+\oplus S^-$ into the two non-isomorphic
  simple $\Cl^{\bar 0}(\frp)$-modules $S^+$ and $S^-$.
    The dual of $S$ is identified with  $\bigwedge U^*$.
  The modules $S$, $S^+$ and $S^-$ are finite dimensional representations of $\widetilde K$ and so are  
their duals and  also $\Cl(\frp)=\bigwedge \frp$, $\Cl(\frp)^{\bar 0}$ and $\Cl(\frp)^{\bar 1}$.

\begin{prop}
As virtual $\widetilde K$ representations, 
\begin{equation}\label{LpSS}
\textstyle \bigwedge^{\bar{\bullet}}\  \frp:=  \bigwedge^{\bar{1}}\frp -  \bigwedge^{\bar{0}} \frp 
 \simeq  (S^+ -S^-)^* \otimes  (S^+ -S^-) \end{equation}
\end{prop}

\proof As  $\widetilde K$-modules :
\[ \textstyle \bigwedge \frp \simeq \Cl(\frp)\simeq \End(S)\simeq S^*\otimes S
 = (S^+\oplus S^-)^*\otimes
(S^+ \oplus S^-).\]
This can be obtained also using :
\[ \textstyle \bigwedge \frp = \bigwedge (U^*\oplus U)\simeq \left( \bigwedge U^*\right)
\otimes \left( \bigwedge U\right) \simeq S^* \otimes S .\]
Writing 
\[(S^+\oplus S^-)^*\otimes (S^+ \oplus S^-)= ((S^+)^*\otimes S^+)\oplus ((S^-)^*\otimes S^+)
\oplus ((S^+)^*\otimes S^-)\oplus ((S^-)^*\otimes S^-)\]
and  identifying the even and odd parts of these $\bbZ/2\bbZ$-graded modules, we get
 \[ \textstyle \bigwedge^{\bar{0}}\frp \simeq ((S^+)^*\otimes S^+)\oplus ((S^-)^*\otimes S^-),\]
 \[ \textstyle \bigwedge^{\bar{1}} \frp \simeq ((S^+)^*\otimes S^-)\oplus ((S^-)^*\otimes S^+).\]
So, as virtual $\widetilde K$ representations,
\[ \textstyle \bigwedge^{\bar{\bullet}}  \frp= (S^+ - S^-)^* \otimes (S^+- S^-).\]
\qed 

\bigskip

\subsection{Dirac cohomology of Harish-Chandra modules}\label{DirCoh}

Let us now  introduce the Dirac operator $D$ :  

\begin{defi}
If $(Y_i)_i$ is a basis of $\frp$ and $(Z_i)_i$ 
is the dual  basis with respect to $B$, then
\[
D=D(\frg,K)=\sum_i Y_i\otimes Z_i\in U(\frg)\otimes \Cl(\frp)
\]
 is independent of the choice of basis  $(Y_i)_i$ and $K$-invariant for the adjoint action 
 on both factors. The Dirac operator  $D$ (for the pair $(\frg,K)$) is an element of $\caA^K$ (see \cite{HPbook}).
\end{defi}

The most important property of $D$ is the formula
\begin{equation}
\label{Dsquared}
D^2=-\mathrm{Cas}_\frg\otimes 1+\Delta(\mathrm{Cas}_{\frk})+(\|\rho_\frk\|^2-\|\rho_\frg\|^2) 
1\otimes 1
\end{equation}
due to Parthasarathy \cite{Par} (see also \cite{HPbook}). Here $\mathrm{Cas}_\frg$ (respectively 
$\mathrm{Cas}_{\frk}$)  denotes the Casimir element of $U(\frg)$. The constant
$(\|\rho_\frk\|^2 -\|\rho_\frg\|^2) $ is explained below.
This formula has several important consequences for Harish-Chandra modules. 
 To state them, we need more notation.
Let us fix a maximal torus $T$ in $K$, with  Lie algebra $\frt_0$.
Let $\fra_0$ denotes the centralizer of $\frt_0$ in $\frp_0$. Then 
\[ \frh_0:= \frt_0\oplus \fra_0 \] 
 is a fundamental Cartan subalgebra of $\frg_0$, and the above decomposition also gives an imbedding
$\frt^* \rightarrow \frh^*$. Let $R=R(\frg,\frh)$ denotes the root system of $\frh$ in $\frg$, 
$W=W(\frg,\frh)$ its Weyl group. Let us also choose a positive root system 
$R^+$ in $R$. As usual, $\rho$ denotes the half-sum of positive roots, an element in $\frh^ *$.
Similarly, we introduce the root system  $R_\frk=R(\frk,\frt)$, its Weyl group $W_\frk$, 
a positive  root system  $R_\frk^+$, compatible with $R^+$, and  half-sum of positive roots $\rho_\frk$.

The bilinear form $B$ on $\frg$ restricts to a non degenerate symmetric  bilinear form
on $\frh$, which is definite positive on the real form $i\frt_0\oplus \fra_0$.
We denote by $\bilo$ the induced form on $i\frt_0^*\oplus \fra_0$  and in the same way its extension to 
$\frh^*$. The norm appearing in (\ref{Dsquared}) is defined for any $\lambda \in \frh^*$
by $\|\lambda\|^2=\bil{\lambda}{\lambda}$. Notice that this is clearly an abuse of notation when 
$\Lambda$ is not in $i\frt_0^*\oplus \fra_0$.

Recall the Harish-Chandra algebra isomorphism 
\begin{equation}\label{HCiso}  \gamma_\frg :\,  \frZ(\frg) \simeq S(\frh)^ W  \end{equation}
between the center $\frZ(\frg)$ of the envelopping algebra $U(\frg)$ and the $W$-invariants in the 
symmetric algebra $S(\frh)$ on $\frh$. Accordingly, a character $\chi$ of $\frZ(\frh)$ is given by 
an element of $\frh^ *$ (or rather its Weyl group orbit). If $\lambda \in \frh^ *$, we denote by 
$\chi_\lambda$ the corresponding character of $\frZ(\frg)$. 
Let $X$ be  a Harish-Chandra module. We say that  $X$ 
has infinitesimal character $\lambda$ if any $z\in \frZ(\frg)$ acts on $X$ by the scalar operator 
$\chi_\lambda(z)\Id_X$. 

Let $\scrM=\scrM(\frg,K)$ be the category of Harish-Chandra modules for the pair $(\frg,K)$ 
(see \cite{KV} for details). If   $X\in \scrM$, then   $\caA=U(\frg)\otimes \Cl(\frp)$ acts 
on $X\otimes S$.   Then $X$ decomposes as the direct sum of its $K$-isotypic components, these 
being finite dimensional if $X$ is admissible.  Accordingly, $X\otimes S$ decomposes 
as the direct sum of its  (finite dimensional if $X$ is admissible) $\widetilde K$-isotypic components.
Let $(\gamma,F_\gamma)$ be an irreducible representation of $\widetilde K$ 
with  highest weight  $\tau=\tau_\gamma\in\frt^*$. We denote the corresponding $\widetilde K$-isotypic
component of $X\otimes S$ by $(X\otimes S)(\gamma)$. Assume $X$ is admissible and has infinitesimal character 
$\Lambda\in\frh^*$.  Then $D^2$ acts on $(X\otimes S)(\gamma)$ by the 
scalar
\begin{equation} \label{scalar}
-\|\Lambda\|^2+\|\tau+\rho_\frk\|^2.
\end{equation}
In particular, we see that in that case  $D^2$ acts semi-simply on $X\otimes S$, and 
 that  the kernel of $D^2$ on $X\otimes S$
is a (finite) direct sum of full $\widetilde K$-isotypic components of $X\otimes S$ : these 
are exactly those $(X\otimes S) (\gamma)$ for which
\begin{equation}
\label{HDKtype}
\|\tau+\rho_\frk\|^2=\|\Lambda\|^2.
\end{equation}

Another important fact is that the action of $D$ preserves $\widetilde K$-isotypic 
components of $X\otimes S$. If $X$ is unitary (resp. finite dimensional), one can put an 
positive definite hermitian form on   $X\otimes S$ and one sees that $D$ is symmetric
(resp. skew-symmetric).

Let us now review Vogan's definition of Dirac cohomology.

\begin{defi} Let $X \in  \scrM$. The Dirac operator $D$ acts on $X\otimes S$ with kernel $\ker D$
 and image $\im D$.
Vogan's Dirac cohomology of $X$ is   the quotient
\[
\Hdir(X) = \ker D/(\ker D\cap \im D ).
\]
\end{defi}

\bigskip 

Since $D \in \caA^{K}$, $\widetilde K$ acts on  $\ker D$,  $\im D$ and $\Hdir(X)$.
Also, assume that $X$ is admissible and  has infinitesimal character $\Lambda \in \frh^ *$.
Then, since $\ker D \subset \ker D^ 2$ and since we have seen that the latter is the sum the full 
$\widetilde K$-isotypic components of $X\otimes S$ satisfying  (\ref{HDKtype}) 
(these are obviously in finite number), we see that 
$\Hdir(X)$ is a finite dimensional representation of $\widetilde K$.
This is particularly helpful if $X$ is unitary, admissible and 
 has infinitesimal character $\Lambda \in \frh^ *$. Then 
it follows that $D$ acts semisimply on $X\otimes S$ and so 
\begin{equation}
\label{HDunit2}
\ker D^2=\ker D = \Hdir(X).
\end{equation}
In this case, the Dirac cohomology of $X$ is a sum the full isotypic components $(X\otimes S)(\gamma)$
such that (\ref{HDKtype}) holds. 
For general $X$, (\ref{HDunit2}) does not hold, but note that $D$ is always a differential
 on $\ker D^2$, and  $\Hdir(X)$ is the usual cohomology of this differential.

Let us state  the main result of \cite{HP}, which gives a strong condition on the 
infinitesimal character of an admissible   Harish-Chandra module $X$ with non zero Dirac cohomology.

\begin{prop}
\label{mainHP} Let $X \in \scrM$ be an admissible  Harish-Chandra module with  infinitesimal character
 $\Lambda\in\frh^*$.  Assume that $(\gamma,F_\gamma)$ is  an irreducible representation of 
$\widetilde K$ with highest weight  $\tau=\tau_\gamma\in\frt^*$ such that    $(X\otimes S)(\gamma) $  
contributes to  $\Hdir(X)$. Then
\begin{equation}
\label{eqHP}
\Lambda=\tau+\rho_\frk\quad\text{up to conjugacy by the Weyl group } W.
\end{equation}
\end{prop}

Thus for unitary $X$, (\ref{HDKtype}) is equivalent to the stronger condition (\ref{eqHP}), provided 
that $\gamma$ appears in $X\otimes S$.

\subsection{Dirac index}

The Dirac index of Harish-Chandra modules 
 is a refinement of Dirac cohomology. It uses the decomposition $ S=S^+ \oplus S^- $
 of the spinor module as a $\Cl(\frp)^{\bar 0}$-module (and thus also as  a representation of  
 $\widetilde K$). Since $D$ is an odd element in $U(\frg)\otimes \Cl(\frp)$ its action
on $X\otimes S$, for any Harish-Chandra module $X$  exchanges $X\otimes S^+$
and $X\otimes S^-$~: 
\[ D: X\otimes S^+  \longleftrightarrow  X\otimes S^-   \] 
Accordingly, the Dirac cohomology of $X$ decomposes as 
\[\Hdir(X)= \Hdir(X)^+\oplus \Hdir(X)^-.\]

The {\bf index of the Dirac operator} acting on $X\otimes S$ is 
the virtual representation
\[ \Idir(X)=\Hdir(X)^+-  \Hdir(X)^- \]
of $\widetilde K$. 
The following proposition is interpreted as an Euler-Poincaré principle.

\begin{prop}\label{IHD}
 Let $X$ be an admissible Harish-Chandra module with infinitesimal character. Then 
\[\Idir(X)=X\otimes S^+ -  X\otimes S^-  \]
as virtual  $\widetilde K$-representations.
\end{prop}

\proof A virtual  $\widetilde K$-representation is by definition an element of 
the Grothendieck group $\scrR( \widetilde K)$ of the category of finite dimensional representations of  
$\widetilde K$. This is the  free $\bbZ$-module generated by equivalence classes of 
irreducible representations, {\sl i.e.} one can write 
\[ \scrR( \widetilde K)=\bigoplus_{\gamma\in (\widetilde K)\, \hat{}} \bbZ. \]
The right-hand side of the equation in the proposition cannot {\sl a priori} be interpreted as an element of 
$\scrR( \widetilde K)$,  but only  of the larger group $\prod_{\gamma\in (\widetilde K)\, \hat{} } \bbZ$.

Let us now prove the equality. We have seen that $D^2$ acts semisimply on $X\otimes S$. 
Furthermore, each eigenspace of 
$D^2$ in $X\otimes S$ is a sum of full $\widetilde K$-isotypic components and that these
are preserved by the action of $D$. Each of these  $\widetilde K$-isotypic components 
$(X\otimes S)(\gamma)$ decomposes also as 
\[ (X\otimes S)(\gamma)=(X\otimes S)(\gamma)^+\oplus (X\otimes S)(\gamma)^-  \]
where $(X\otimes S)(\gamma)^\pm:=(X\otimes S)(\gamma) \cap (X\otimes S^\pm)$.

For $\widetilde K$-isotypic components corresponding to a non-zero eigenvalue of $D^2$, we thus get 
that $D$ is a bijective intertwining operator (for the $\widetilde K$-action) between
$(X\otimes S)(\gamma)^+$ and $(X\otimes S)(\gamma)^-$. Thus the contribution of these
 $\widetilde K$-isotypic components to $ X\otimes S^+ -  X\otimes S^- $ is zero. So only $\ker D^2$
will contribute, {i.e.}
\[X\otimes S^+ -  X\otimes S^-= (\ker D^2 \cap (X\otimes S^+))-(\ker D^2 \cap (X\otimes S^-)). \]
Let us write 
\[ \ker D^2 \cap (X\otimes S^\pm) =\ker D \cap (X\otimes S^\pm)\oplus W^\pm  \]
for some $\widetilde K$-invariant complementary subspaces $W^\pm$.
Then, as above, $D$ is  bijective intertwining operator for the $\widetilde K$-action
between $W^\pm$ and  $D(W^\pm) \subset \ker D^2 \cap (X\otimes S^\mp)$. So these contributions
 also cancel, and what remains is exactly the virtual $\widetilde K$-representation
 $\Hdir(X)^+- \Hdir(X)^-$. \qed

\bigskip

\subsection{Dirac index of limits of discrete series}\label{LDS}

We assume that $G$ has discrete series, so that $G$ and $K$ have same rank, and 
$T$ is a compact Cartan subgroup. Harish-Chandra modules of limits of discrete series of $G$ are obtained as cohomologically induced
$A_\frb(\bbC_\lambda)$-modules (see \cite{KV}), where 
$\frb=\frt\oplus \fru$ is a Borel subalgebra containing $\frt$ with nilpotent radical $\fru$, 
and $\bbC_\lambda$ is the one-dimensional representation of $T$ with weight
$\lambda \in i\frt_0^*$. Some positivity conditions on $\lambda$ are required, that we now describe.
 The Borel subalgebra $\frb$ determines a set of positive roots 
$R^+_\frb$ of $R=R(\frt,\frg)$ (the roots of $\frt$ in $\fru$). Let us denote
by $\rho(\frb)$, (resp.  $\rho_c(\frb)$,  resp. $\rho_n(\frb)$) the half-sum of (resp. compact , 
resp.non-compact)  roots in $R^+_\frb$.  The positivity condition on $\lambda$ is that 
 \begin{equation} \label{posi} \bil{\lambda+\rho(\frb)}{\alpha}\geq 0.\end{equation}
 Then, $A_\frb(\bbC_\lambda)$ is a discrete series modules  if the inequalities in (\ref{posi})
 are strict, and otherwise a limit of discrete series ($\chi=\lambda+\rho(\frb)$ is the infinitesimal 
 character of    $A_\frb(\bbC_\lambda)$), or $0$ (but we are not interested in this case). 
 The lowest $K$-type of  $A_\frb(\bbC_\lambda)$
 has multiplicity one and highest weight $\Lambda=\lambda+2\rho_n(\frb)=\chi+\rho_n(\frb)-\rho_c(\frb)$, 
 and all other $K$-types have highest weights  of the form $\Lambda+ \sum_j\beta_j$ for some positive 
 roots  $\beta_j$.

 In \cite{HKP} these facts are used together with Prop \ref{mainHP}  to show that  the Dirac cohomology   
  of  $A_\frb(\bbC_\lambda)$ consists in the multiplicity-one $\widetilde K$-type $F_\mu$ with 
   highest weight 
   \[\mu=\Lambda-\rho_n(\frb)=\lambda+\rho_n(\frb)=\chi-\rho_c(\frb).\]
   (In \cite{HKP}, it is assumed that $A_\frb(\bbC_\lambda)$ is a discrete series, but 
   inspection of the proof easily shows that it works also for limits of discrete series.)
Thus the Dirac index of    $A_\frb(\bbC_\lambda)$ is $\pm F_\mu$. 
To determine the sign, recall that the spinor module $S$, as a $\widetilde K$-representation, 
doesn't depend on any choice, nor does the set $\{S^+, S^-\}$, in particular not on the way we 
realized this module, but the  distinction between   $S^+$ and $S^-$ does (resulting  
on a sign change in the Dirac index). So suppose we have fixed once for all
a Borel subalgebra $\frb_1=\frt\oplus \fru_1$ and choose $\fru_1\cap \frp$ as the isotropic subspace
$U$ of $\frp$ used to construct the spinor modules in section \ref{gK}.
Then $S^+$ is the $\widetilde K$-representation containing the weight $-\rho_n(\frb_1)$ and 
$S^-$ is the one  containing the weight $\rho_n(\frb_1)$. With this choice, it is easy 
to see that $\Idir(A_\frb(\bbC_\lambda))=\mathrm{sgn}(w) F_\mu$, where $w \in W$ is the Weyl group element 
sending $\frb_1$ to $\frb$ and $\mathrm{sgn}$ is the sign character.

Let us now determine the Dirac index of the virtual modules $X$ which are the linear combinations 
of limits of discrete series with the same infinitesimal character, whose characters are the 
supertempered distributions constructed by Harish-Chandra (see \cite{HC} and \cite{Bouaziz}, \S 7).
Namely take an integral but non necessarilly regular weight $\chi$ in $i\frt_0^*$
and consider the  limits of discrete   series $A_\frb(\bbC_\lambda)$ as above 
with $\chi=\lambda+\rho(\frb)$ satisfying (\ref{posi}). Notice that $\chi$ being fixed, 
$\frb$ determines $A_\frb(\bbC_\lambda)$ and that  Borel subalgebras $\frb$ 
which occurs are the one such that the corresponding Weyl chamber $C_\frb$  has $\chi$ in its closure.
So we can forget the $\lambda$ in the notation, and the set of  limits of discrete   series 
we are considering  is $\caL=\{ A_\frb\, \vert \chi \in \overline{C_\frb}   \}$. 
We should also take care of the fact that the $A_\frb(\bbC_\lambda)$ could be $0$, but 
the important remark here is that if one of them in the set above is non zero, then all are 
($\chi $ is  not on a wall in the kernel of  a simple compact root).
Choose one of them as a base point, say $A_{\frb_2}$. Then the linear combination introduced by Harish- 
Chandra is 
\[ X_{\chi, \frb_2}=\frac{1}{\vert W_\chi \vert } \sum_{A_\frb \in \caL } \epsilon(\frb)  A_\frb   \]
where $\epsilon(\frb)=  \mathrm{sgn}(w)$, and  $w \in W$ is the Weyl group element sending
 $\frb_2$ to $\frb$ and $W_\chi=\{w\in W|\, w\cdot \chi=\chi\}$. This construction is made because all the  the $\epsilon(\frb)  A_\frb$ 
 have the same character formula   on $T$ and this is also the formula for 
 the character of $X_{\chi, \frb_2}$.   Notice that a different  choice of $\frb_2$ as  base 
 point would result in at most a sign change.

The fact that the $A_\frb$ are non-zero implies that the corresponding Weyl chambers $C_\frb$ are all 
included in a single Weyl chamber for $R(\frk,\frt)$. In particular, the various
$\rho_c(\frb)$ are all equal (let say to $\rho_c$). This shows that the Dirac cohomology of all the 
$  A_\frb \in \caL$ is the same, namely the multiplicity-one $\widetilde K$-type $F_\mu$ with 
   highest weight $\mu=\chi-\rho_c$. Taking signs into account, we see that all the 
 $\epsilon(\frb)  A_\frb $ have same Dirac index $\mathrm{sgn}(w) F_\mu$
 where $w$ is the Weyl group element sending $\frb_1$ to $\frb_2$. This is thus also the Dirac index of 
 $X_{\chi,\frb_2}$. 

For discrete series, $\chi$ is regular, so $\caL$ contains only one element $A_\frb(\lambda)$
with $\lambda=\chi-\rho(\frb)$. 

\begin{prop}\label{OrtInd}
If $\Idir(X_{\chi,\frb})=  \pm I(X_{\chi',\frb'})$ then $\chi$ and $\chi'$
are conjugate by a element of $W(\frk,\frt)$, and thus by an element in $K$.
Therefore  $X_{\chi,\frb}=  X_{\chi',\frb'}$. 
\end{prop}

\proof    Say that  $\Idir(X_{\chi,\frb})= \pm F_\mu=  \pm I(X_{\chi',\frb'})$. 
For any   $g\in G$ normalizing $T$, $X_{g\cdot \chi,  g \cdot \frb}= X_{\chi,\frb}$, 
so we can assume that 
$C_\frb$ and $C_{\frb'}$ are in the same Weyl chamber for $R(\frk,\frt)$.
We have then $\mu=\chi-\rho_c(\frb)=\chi'-\rho_c(\frb')$ and  $\rho_c(\frb)=\rho_c(\frb')$ 
so that $\chi=\chi'$.  \qed

\section{Euler-Poincaré pairing and Dirac pairing}\label{EPdir}

For two finite length  Harish-Chandra modules $X$ and $Y$, one can define their 
Euler-Poincaré pairing as the alternate sum of dimensions of $\Ext$ functors as 
in (\ref{EPp}).
Recall that $\scrM$ has finite cohomological dimension, so this sum has finite support.
 More precisely, an  explicit projective resolution of $X$ in $\caM(\frg,K)$
is given by  
\[\cdots  \longrightarrow P_{i+1} \longrightarrow P_{i}  \longrightarrow \cdots \longrightarrow P_{0}
 \longrightarrow X \longrightarrow 0\]
with $P_i =(U(\frg)\otimes_{U(\frk)}\bigwedge^ i \frp)\otimes X$ (\cite{KV}, \S II.7). 
Set 
\[C^i=\Hom_\scrM(P_i, X)\simeq   \Hom_{K}( \textstyle {\bigwedge^i }\frp\otimes X, Y).\]
Thus   $\Ext^i(X,Y)$ is given by the $i$-th cohomology group of the complex ${\bf C}=(C^i)_i$, with differential $d^i$
given explicitely in {\sl loc. cit.}

If $\gamma, \sigma$ are virtual finite dimensional representations of $\widetilde K$,
and $\chi_\gamma$, $\chi_\sigma$ are  their characters, we denote by  
\[\left[\gamma, \sigma \right]_{\widetilde K}=\int_{\widetilde K} \overline{\chi_\gamma(k)} \chi_\sigma(k)\; dk \]
 the usual (hermitian) pairing between these virtual representations ($dk$ is the normalized 
 Haar measure on $\widetilde K$).
Since  virtual finite dimensional  representations of $K$ are also 
 virtual  representations of $\widetilde K$, we use also the notation for their pairing.
Notice that when $\gamma, \sigma$ are actual finite dimensional representations, then 
\[\left[\gamma, \sigma \right]_{\widetilde K}=\dim \Hom_{\widetilde K}(\gamma, \sigma).\]
We can now compute, assuming that $X$ and $Y$ have infinitesimal character \footnote{A technical remark is in order here : 
recall that a Harish-Chandra  module with infinitesimal character is of finite length if and only if it is admissible. } : 

 \begin{align*}
 \EP(X,Y)&=\sum_i (-1)^i \; \dim \Ext^i(X,Y) \\
&=  \sum_i (-1)^i \;  \dim  \Hom_{K}\left(\textstyle {\bigwedge^i} \frp\otimes X, Y\right)  \quad  (\text {Euler-Poincaré principle}) \\
   &=  \sum_i (-1)^i \; \left[ \textstyle{\bigwedge^i} \frp \otimes X, Y\right]_{\widetilde K}
   =  \left[ \sum_i (-1)^i \; \textstyle{\bigwedge^i} \frp\otimes X, Y\right]_{\widetilde K}
=  \left[ \textstyle{\bigwedge^{\bar{\bullet}}} \frp\otimes X, Y\right]_{\widetilde K} \\
   &=  \left[(S^ +-S^-)^*\otimes (S^ +-S^-) \otimes X, Y\right]_{\widetilde K} \qquad 
(\text{see } (\ref{LpSS}))\\
  &=   \left[ (S^ +-S^-)  \otimes X, (S^ +-S^-)\otimes Y\right]_{\widetilde K}  \\
  &=   \left[ \Hdir(X)^+ -  \Hdir(X)^- ,
\Hdir(Y)^+-  \Hdir(Y)^- \right]_{\widetilde K} \qquad  (\text{Prop. } \ref{IHD}.)\\
 &=   \left[ \Idir(X), \Idir(Y) \right]_{\widetilde K}.
\end{align*}

The attentive reader probably noticed a small problem with the above computation, namely  the application of 
the Euler-Poincaré principle is not justified since  the terms $C^i=  \Hom_{K}\left(\textstyle {\bigwedge^i}
 \frp\otimes X, Y\right) $ could be infinite dimensional.  As explained in the introduction, we were not able
 to find a proof of the equality of the two extreme terms so far, so let us state this as a conjecture :

 \begin{conj}\label{Main1} Let $X$ and $Y$ be  two finite length  Harish-Chandra modules with infinitesimal character.
 Then 
 \[\EP(X,Y)= \left[ \Idir(X), \Idir(Y) \right]_{\widetilde K}= \bil{X}{Y}_{\mathrm{Dir}}.  \]
 \end{conj}
The last equality is the  definition of the  Dirac pairing mentioned in the introduction. 
The main content of this conjecture   is that the $\EP$ pairing of two finite length Harish-Chandra 
modules  with infinitesimal character factors through their Dirac indices, and thus through their Dirac cohomology.
In particular,  the results on Dirac cohomology and Dirac index 
recalled in Section \ref{DirCoh} put severe conditions on modules $X$ and $Y$ for their Dirac  pairing to be 
non-zero. For many interesting modules, the Dirac index is explicitly known, and then so is 
 the Dirac pairing between these modules.

\bigskip
\begin{rmk}
In the case where $X=F$ is a finite-dimensional irreducible Harish-Chandra module, $\Ext^i(F,Y)$
is also the $i$-th $(\frg,K)$-cohomology group of $Y\otimes F^*$.
Thus, in that case, 
\[\sum_i (-1)^i \dim H^i(\frg, K; Y\otimes F^*)=  \left[ \Idir(F), \Idir(Y) \right]_{\widetilde K}.\]
In the case where $Y$ is unitary, we have much stronger results : 
the differential $d$ on $\Hom_K(\bigwedge \frp \otimes F, Y)$ is $0$ 
so 
\[\bigoplus_{i}  H^i(\frg, K; Y\otimes F^*)=\bigoplus_{i} \Hom_K(\textstyle \bigwedge^i \frp \otimes F, Y)= 
\Hom_K(\bigwedge \frp \otimes F, Y) \simeq \Hom_{\widetilde K}(F\otimes S, Y\otimes S) \]  
and the only common $\widetilde K$-types between $F\otimes S$ and $Y\otimes S$
have their isotypic components in $\ker (D_{\vert F\otimes S})=\ker (D_{\vert F\otimes S}^2)=\Hdir(F)$
and $\ker (D_{\vert X\otimes S})=\ker (D_{\vert X\otimes S}^2)=\Hdir(X)$ respectively. 
See \cite{Wal}, \S 9.4  and \cite{HPbook},  \S 8.3.4.
\end{rmk}

\section{Labesse index functions and Euler-Poincaré functions} \label{EPLab}

Let us continue the computation of  the  Dirac pairing (and thus conjecturally  of the   $\EP$ pairing)
 for two finite length  Harish-Chandra 
modules $X$ and $Y$ with infinitesimal character. We have 
\[ \bil{X}{Y}_{\mathrm{Dir}} = 
    \left[ \Idir(X), \Idir(Y) \right]_{\widetilde K}  = \sum_{\gamma\in (\widetilde K)\, \hat{}} 
\left[\gamma, \Idir(X) \right]_{\widetilde K}  \times 
\left[ \gamma, \Idir(Y) \right]_{\widetilde K}.\]

Let us now introduce some material from \cite{Lab}. If $\gamma$ is a genuine 
virtual finite dimensional representation of $\widetilde K$, let us denote 
by $I_K(\gamma)$ the virtual finite dimensional representation
\[ \gamma\otimes (S^+  -S^-)^*   \]
of $K$ and by  $\tilde  \chi_\gamma$ of  the character of $I_K(\gamma)$ (a conjugation invariant function on $K$).
If $(\pi,X)$ is an admissible Harish-Chandra module, with $\pi$ denoting the action of $K$ on $X$, 
 the operator 
\[  \pi(\overline{\widetilde{\chi_\gamma}})=\int_K \pi(k) \,  \overline{\tilde \chi_\gamma (k)} \; dk  \]
is of finite rank, and so is a trace operator. Let us denote its trace by 
$I(X,\gamma)$. If furthermore  $X$ has infinitesimal character,  we have
\[I(X,\gamma)=\left[  \gamma\otimes (S^+-S^-)^*  , X \right]_{\widetilde K}=  
\left[  \gamma , X \otimes (S^+-S^-)\right]_{\widetilde K}= \left[  \gamma , \Idir(X)\right]_{\widetilde K}.\]

Labesse main result is the existence   for all 
$\gamma\in (\widetilde K)\, \hat{}$  of a 
smooth, compactly supported and bi-$K$-finite function $f_{ \gamma}$
which satisfies, for all finite length  Harish-Chandra module $X$
\[ \Theta_X(f_\gamma)= I(X,\gamma),   \]
where $\Theta_X$ denotes the distribution-character of $X$. The main ingredient in Labesse's construction is 
the Paley-Wiener theorem of Arthur (\cite{Ar}).

From this, we get 
\[  \bil{X}{Y}_{\mathrm{Dir}} = \sum_{\gamma\in (\widetilde K)\, \hat{}} I(X,\gamma)\times I(Y,\gamma)
= \sum_{\gamma\in (\widetilde K)\, \hat{}}   \Theta_X(f_\gamma) 
 \Theta_Y(f_\gamma). \]

Let us call
\[ f_Y=  \sum_{\gamma\in (\widetilde K)\, \hat{}}   \Theta_Y(f_\gamma) 
\;  f_\gamma   \]
an index (or Euler-Poincaré, if we believe in Conjecture \ref{Main1})  function for $Y$. 
Note that the sum above has finite support, so that 
$f_Y$ is smooth, compactly supported and bi-$K$-finite. Furthermore, Labesse shows
(\cite{Lab}, Prop. 7) that the functions $f_\gamma$ are "very cuspidal", {\sl i.e.}
that their constant terms for all proper parabolic subgroups of $G$ vanish.
Thus the same property holds for $f_Y$. This implies the vanishing of the orbital integrals 
$\Phi(f_Y,x)$ on  regular non-elliptic element $x$ in $G$. We have obtained : 

\begin{thm}
For any finite length Harish-Chandra module $Y$ with infinitesimal character, there exists a 
smooth, compactly supported,   bi-$K$-finite  and very cuspidal function
$f_Y$ on $G$ such that for any finite length Harish-Chandra module $X$ with infinitesimal character,
\[  \bil{X}{Y}_{\mathrm{Dir}} = \Theta_X(f_Y). \]
\end{thm}

 \section{Integral orbital of $f_Y$ as the character of $Y$ on elliptic elements}\label{OrbInt}

The goal of this section is to show that the value of the orbital integral $\Phi(f_Y,x)$ of $f_Y$ at 
an elliptic regular element $x$ coincide with the value of the character $\theta_Y$ of $Y$ at $x^{-1}$.
 The character $\Theta_Y$ of $Y$ is a distribution on $G$,  but recall that 
according to Harish-Chandra regularity theorem, there is an analytic,  conjugation invariant function
that we will denote  by $\theta_Y$ on the set $G_{\mathrm{reg}}$ such that for
 all $f\in\caC^\infty_c(G)$, 
\[\Theta_Y(f)=\int_G  \theta_Y(x) \; f(x)\; dx.    \]
Using Weyl integration formula, this could be written as 
\[\Theta_Y(f)=\sum_{[H]} \frac{1} {\vert  W(G,H)  \vert } \int_{H}   \vert D_G(h)\vert \, 
  \theta_Y(h) \;     \Phi(f, h)\; dh    \]
where the first sum is on a system of representative of conjugacy classes of Cartan subgroups $H$ of $G$, 
 $W(G,H)=N_G(H)/H$ is the real Weyl group of $H$ and $\vert D_G\vert$ is the usual Jacobian.
The assertion is thus that :
\begin{equation}\label{assertion} \theta_Y(x^{-1})=\Phi(f_Y,x), \qquad 
(x\in G_{\mathrm{ell}}) .\end{equation}
\medskip 

The characterization of orbital integrals due to A. Bouaziz \cite{Bouaziz} shows that there indeed 
exists a function $\psi_Y$ in $\caC^\infty_c(G)$ such that 
 \[  \Phi(\psi_Y,x)=\begin{cases}  \theta_Y(x^{-1})  \text{ if } x\in G_{\mathrm{ell}} \\
 0 \text{ if } x\in G_{\mathrm{reg}}\setminus G_{\mathrm{ell}}  \end{cases}.\]

Let $\caF$ be a family of elements $X$ in $\scrR_\bbC$ such that the space generated by the $\Theta_X$ 
is dense  in the space of invariant distributions $\caD'(G)^G$ on $G$. For instance, $\caF$ could be the 
set of  characters of all  irreducible tempered representations, but we  will rather  take $\caF$ to be the  family of 
 virtual representations    with characters  $\Theta_{h^*}$   defined in  \cite{Bouaziz}, \S 7. 
 The density of this family    of invariant distributions is a consequence of the inversion
  formula of orbital 
  integrals \cite{Bouaziz2}.   Elements in  $\caF$  are generically irreducible tempered representations,
   but in general, they are  linear combinations of some of these with same infinitesimal character.
   The distributions $\Theta_{h^*}$ are supertempered in the sense of Harish-Chandra \cite{HC}.
In any case, by density of the family $\caF$,  to prove (\ref{assertion}), it is  enough to show that for
 all $X\in \caF$,  we have : 
\begin{equation}\label{reduc} \Theta_X(f_Y)=\Theta_X(\psi_Y).\end{equation}
If $X$ is a linear combination of parabolically induced representations,  then both side are $0$ since 
$\Theta_X$ vanishes on elliptic elements. Thus,  it is sufficient to prove (\ref{reduc}) for $X$
corresponding to the  $\Theta_{h^*}$ of \cite{Bouaziz} attached to the fundamental Cartan subgroup $T$.
For simplicity, we assume now that $T$ is compact, {\sl i.e.} that $G$ and $K$ have same rank. Then  
 $X$ is either a  discrete series or a linear combination of limits of discrete series (with same infinitesimal character) described in section \ref{LDS}.

Assume that (\ref{assertion}) is established for all such $X$.
The left-hand side of (\ref{reduc}) then  also  equals 
\[\Theta_Y(f_X)=\frac{1} {\vert  W(G,T)  \vert } \int_{T}   \vert D_G(t)\vert \,   \theta_Y(t) \;     \Phi(f_X,t ) \; dt\]
\[=\frac{1} {\vert  W(G,T)  \vert } \int_{T}   \vert D_G(t)\vert\,   \theta_Y(t) \;     \theta_X(t^{-1} ) \; dt  \]
and  equals the right-hand side, by using the definition of $\psi_Y$ and  the Weyl integration formula again.
By definition of the measure $dc$ on the set $C^{\mathrm{ell}}$  of regular semisimple
elliptic conjugacy classes in $G$ and of the elliptic pairing recalled in the introduction, we have   also 
\[\Theta_Y(f_X)= \int_{  C^{\mathrm{ell}}  } \theta_Y(c) \;     \theta_X(c^{-1}) \; dc=\bil{X}{Y}_{\mathrm{ell}}.\] 

Thus we have reduced the proof of (\ref{assertion}) for all $Y$ to the case when 
$Y$ is either a  discrete series or a linear combination of limits of discrete series as described above.
In turns, it is enough to show (\ref{reduc}) when both $X$ and $Y$ are of this kind.
In case $Y$ corresponds to a parameter $h^*$ of \cite{Bouaziz2}, $\psi_Y$ is exactly 
the function denoted $\psi_{h^*}$ there, and in particular : 
\[\Theta_X(\psi_Y)=1 \text{ if } X=Y, \quad 0 \text{ otherwise}. \]

With the notation of section \ref{LDS}, we can take $X=X_{\chi,\frb}$ and $Y=X_{\chi',\frb'}$.
Then $f_X=\sum_{\gamma\in (\widetilde K)\, \hat{}} \Theta_X(f_\gamma)f_\gamma=   \sum_{\gamma\in (\widetilde K)\, \hat{}}    \left[\gamma,  I(X,\gamma) \right]_{\widetilde K} \, f_\gamma$
and by the results in section \ref{LDS}, we get $f_X=\epsilon(\frb) f_{\chi-\rho_c(\frb)}$.
Similarly we get $f_Y=\epsilon(\frb') f_{\chi'-\rho_c(\frb')}$ and 
 \[ \Theta_X(f_Y)= \left[ \Idir(X), \Idir(Y) \right]_K = \left[ \epsilon(\frb) \, F_{\chi-\rho(\frb)},
 \epsilon(\frb')\, F_{\chi'-\rho(\frb')} \right]_K .\] 
 By Proposition \ref{OrtInd}, we see that this is $0$ if $X\neq Y$ and $1$ if $X=Y$. 
 This finishes the proof of  (\ref{reduc}) in  the case under consideration. 
We have proved : 

\bigskip

\begin{thm}
Let $X$ and $Y$ be finite length Harish-Chandra modules with infinitesimal character in $\scrM$ and let $f_X$, $f_Y$ be 
the Euler-Poincaré functions for $X$ and $Y$ respectively,  constructed in section \ref{EPLab}. 
Then, the orbital integral $\Phi(f_X,x)$ at a regular  element $x$ of $G$ is $0$ if 
$x$ is not elliptic, and equals $ \theta_X(x^{-1})$
if $x$ is elliptic.
Furthermore :
\[ \bil{X}{Y}_{\mathrm{Dir}}=\Theta_X(f_Y)=\Theta_Y(f_X)= \int_{  C^{\mathrm{ell}}  } \theta_Y(c) \;     \theta_X(c^{-1}) \; dc= \bil{X}{Y}_{\mathrm{ell}}. \]
\end{thm}

\section{About Conjecture \ref{Main1}}\label{Proof}

As explained in the introduction, the goal of this section is to provide evidence for Conjecture \ref{Main1}, as 
well as a possible direction of investigation for a proof.
In what follows, $X$ and $Y$ are finite length Harish-Chandra modules with infinitesimal characters.
We may assume that the infinitesimal characters of $X$ and $Y$ are the same, since otherwise, both side
of the identity we aim to prove are $0$.
\bigskip 

Let $\scrC=\Hom_{\widetilde K}(X\otimes S, Y\otimes S)$. Then $\scrC=\scrC^{\bar 0}\oplus \scrC^{\bar 1}$, with 
\[\scrC^{\bar 0}= \Hom_{\widetilde K}(X\otimes S^+, Y\otimes S^+)\oplus 
\Hom_{\widetilde K}(X\otimes S^-, Y\otimes S^-),   \]
\[\scrC^{\bar 1}= \Hom_{\widetilde K}(X\otimes S^+, Y\otimes S^-)\oplus 
\Hom_{\widetilde K}(X\otimes S^-, Y\otimes S^+).  \]

Let us consider the following various actions of the Dirac operator : 
\[D_X^{+-} :   X\otimes S^+\longrightarrow X\otimes S^-, \qquad
D_X^{-+} :   X\otimes S^-\longrightarrow X\otimes S^+, \]
\[D_Y^{+-} :   Y\otimes S^+\longrightarrow Y\otimes S^-, \qquad
D_Y^{-+} :   Y\otimes S^-\longrightarrow Y\otimes S^+, \]

\medskip

For  $\phi^{++}\in \Hom_{\widetilde K}(X\otimes S^+, Y\otimes S^+)$,  and  
$\phi^{--}\in \Hom_{\widetilde K}(X\otimes S^-, Y\otimes S^-)$ set
\[ S\phi^{++}=-\phi^{++}\circ D_X^{-+}+D_Y^{+-}\circ \phi^{++}.\]
\[ S\phi^{--}=-\phi^{--}\circ D_X^{+-}-D_Y^{-+}\circ \phi^{--}.\]
This defines a linear map  $S: \scrC^{\bar 0} \rightarrow \scrC^{\bar 1}$.

\medskip

For  $\psi^{+-}\in \Hom_{\widetilde K}(X\otimes S^+, Y\otimes S^-)$,  and 
 $\psi^{-+}\in \Hom_{\widetilde K}(X\otimes S^-, Y\otimes S^+)$, set
\[ T\psi^{+-}=-\psi^{+-}\circ D_X^{-+}+D_Y^{-+}\circ \psi^{+-}.\]
\[ T\psi^{-+}=-\psi^{-+}\circ D_X^{+-}-D_Y^{+-}\circ \psi^{-+}.\]
This defines a linear map  $T : \scrC^{\bar 1} \rightarrow \scrC^{\bar 0}$.

\medskip

Let us take $\phi^{++}+\phi^{--}$ in $\ker S$. Thus 
\begin{equation}\label{kerS} 
\phi^{++}\circ D_X^{-+}+D_Y^{-+}\circ \phi^{--}=0 \text{ and } -\phi^{--}\circ D_X^{+-}+D_Y^{+-}\circ \phi^{++}=0.\end{equation}
 From this we see that 
\[ \phi^{++}(\ker(D_X^{+-}))\subset \ker(D_Y^{+-})  \text{ and } \phi^{--}(\ker(D_X^{-+}))\subset \ker(D_Y^{-+}) ,\]
and  also that  
\[  \phi^{++}(\im(D_X^{-+}))\subset \im(D_Y^{-+}) \text{ and } \phi^{--}(\im(D_X^{+-}))\subset \im(D_Y^{+-}).  \]
Therefore $\phi^{++}$ induces
\[\bar \phi^{++}:\Hdir(X)^+= \frac{\ker (D_X^{+-}) }{ \ker (D_X^{+-})\cap  \im(D_X^{-+})  }
\longrightarrow \Hdir(Y)^+=\frac{\ker (D_Y^{+-}) }{\ker (D_Y^{+-})\cap  \im(D_Y^{-+})} , \]
and 
 $\phi^{--}$ induces
\[\bar \phi^{--}: \Hdir(X)^-=\frac{\ker (D_X^{-+}) }{ \ker (D_X^{-+})\cap  \im(D_X^{+-})}  
\longrightarrow \Hdir(Y)^-=\frac{\ker (D_Y^{-+}) }{ \ker (D_Y^{-+})\cap  \im(D_Y^{+-})} . \]

Let us now show that if $\phi^{++}+\phi^{--}\in \ker S\cap \im T$, then $(\bar \phi^{++}, \bar \phi^{--})=(0,0)$. 
Write 
\[ \phi^{++}=  - \psi^{-+}\circ D_X^{+-}+D_Y^{-+}\circ \psi^{+-}, \quad \phi^{--}= 
- \psi^{+-}\circ D_X^{-+}-D_Y^{+-}\circ \psi^{-+}.\]
This implies that 
\begin{equation}\label{incl}\textstyle
  \phi^{++}(\ker(D_X^{+-}))\subset \im(D_Y^{-+})\cap \ker(D_Y^{-+}),  \;
      \phi^{--}(\ker(D_X^{-+}))\subset \im(D_Y^{+-})\cap \ker(D_Y^{+-})  ,\end{equation}
which proves the assertion. Therefore, 
we have shown that there is a well-defined morphism $\phi^{++}+\phi^{--}\mapsto  \bar \phi^{++}+\bar\phi^{--}$
from $\displaystyle \frac{ \ker S}{ \ker S\cap \im T} $ to 
\[ \Hom_{\widetilde K}  \left( \Hdir(X)^+,\Hdir(Y)^+\right) 
\oplus 
\Hom_{\widetilde K}  \left(  \Hdir(X)^-  ,\Hdir(Y)^-  \right) .   \]

Similarly, there is  a well-defined morphism $\psi^{+-}+\psi^{-+}\mapsto  \bar \psi^{+-}+\bar\psi^{-+}$
from $\displaystyle \frac{ \ker T}{ \ker T\cap \im S}$ to 
\[ \Hom_{\widetilde K}  \left( \Hdir(X)^+,\Hdir(Y)^-\right) 
\oplus 
\Hom_{\widetilde K}  \left(  \Hdir(X)^-  ,\Hdir(Y)^+  \right) .   \]

\bigskip 

\begin{lemma}\label{chiant}  The two  morphisms defined above are isomophisms : 
 
\[ \frac{\ker S}{ \ker S\cap \im T}\simeq \Hom_{\widetilde K}  \left( \Hdir(X)^+,\Hdir(Y)^+\right) 
\oplus 
\Hom_{\widetilde K}  \left(  \Hdir(X)^-  ,\Hdir(Y)^-  \right),  \] 
\[\frac{\ker T}{ \ker T\cap \im S }\simeq  \Hom_{\widetilde K}  \left( \Hdir(X)^+,\Hdir(Y)^-\right) 
\oplus 
\Hom_{\widetilde K}  \left(  \Hdir(X)^-  ,\Hdir(Y)^+  \right) .\] 
\end{lemma}
\bigskip 

\proof 
Let us compute $TS$. If $\phi^{++}\in \Hom_{\widetilde K}(X\otimes S^+, Y\otimes S^+)$, we
get 
\begin{align*}& TS\phi^{++} = T(-\phi^{++}\circ D_X^{-+}+D_Y^{+-}\circ \phi^{++})\\
&=   \phi^{++}\circ D_X^{-+}\circ D_X^ {+-}+D_Y^{+-} \circ \phi^{++}\circ D_X^{-+}
-D_Y^{+-}\circ \phi^{++}\circ  D_X^{-+}+D_Y^{-+}\circ D_Y^{+-}\circ \phi^{++}\\
&= \phi^{++}\circ D_X^{-+}\circ D_X^ {-+}+D_Y^{-+}\circ D_Y^{+-}\circ \phi^{++}, 
\end{align*}
and if $\phi^{--}\in \Hom_{\widetilde K}(X\otimes S^-, Y\otimes S^-)$, 
\begin{align*}& TS\phi^{--} = T(-\phi^{--}\circ D_X^{+-}-D_Y^{-+}\circ \phi^{--})\\
& =   \phi^{--}\circ D_X^{+-}\circ D_X^ {-+}- D_Y^{-+} \circ \phi^{--}\circ D_X^{+-}
       +D_Y^{-+}\circ \phi^{--}\circ  D_X^{+-}+D_Y^{+-}\circ D_Y^{-+}\circ \phi^{--}\\
&= \phi^{--}\circ D_X^{+-}\circ D_X^ {+-}+D_Y^{+-}\circ D_Y^{-+}\circ \phi^{--}.
\end{align*}
Similar computations for $ST$ show that  in fact, for any $\phi \in \scrC^{\bar 0}$, and any $\psi \in  \scrC^{\bar 1} $
\begin{equation} TS\phi=D^2\circ \phi +\phi \circ D^2, \quad  ST\psi=D^2\circ \psi +\psi \circ D^2
\end{equation}
Let us now use the fact that 
\[  \scrC=\Hom_{\widetilde K}(X\otimes S,Y\otimes S)=\prod_{\gamma \in (\widetilde K)\hat{}}
\Hom_{\widetilde K}((X\otimes S)(\gamma),(Y\otimes S)(\gamma)).\]
Since $X$ and $Y$ have same infinitesimal character, by Equation (\ref{scalar}), 
$D^2$ acts on $(X\otimes S)(\gamma)$ and  $(Y\otimes S)(\gamma)$
by the same  scalar, let us call it $\alpha_\gamma$ for short,  so we see that on each $\bbZ/2\bbZ$-graded component of 
$\Hom_{\widetilde K}((X\otimes S)(\gamma),(Y\otimes S)(\gamma))$, $TS$ or $ST$ (depending on the  component)
is multiplication by $\alpha_\gamma^2$.
From this, we deduce two facts that are worth recording.
\begin{lemma} \label{plustard}$\ker (TS)$ (resp. $\ker ST$) is finite dimensional, and 
\[ker (TS)\oplus \im (TS)=\scrC^{\bar 0}, \quad \ker(ST)\oplus \im(ST)=\scrC^{\bar 1}.\]
\end{lemma}
Indeed, we see that $\ker(TS)$ (or $ST$) is 
\[ \prod_{\gamma \in (\widetilde K)\, \hat{}\; \vert \; (X\otimes S)(\gamma)\subset \ker D^2 }
\Hom_{\widetilde K}((X\otimes S)(\gamma),(Y\otimes S)(\gamma)).\] 
Notice that  $\scrC^{\bar 0}$ may not  be written as a direct sum of eigenspaces for $TS$, but is 
isomorphic to a direct product for which each factor is an eigenspace for $TS$. 

\bigskip 

Let us show now that the morphisms in Lemma \ref{chiant} are injective and surjective.
We will do it only for the first one, so suppose  that $\phi=\phi^{++}+\phi^{--}\in \ker S$ is such that 
$(\bar \phi^{++}, \bar\phi^{--})=(0,0)$. 
This is equivalent to (\ref{incl}). We want to find $\psi=\psi^{+-}+\psi^{-+}$ such that 
$\phi=T\psi$, {\sl i.e.}
\begin{equation}\label{WWW} \phi^{++}=  - \psi^{-+}\circ D_X^{+-}+D_Y^{-+}\circ \psi^{+-}, \quad \phi^{--}= 
- \psi^{+-}\circ D_X^{-+}-D_Y^{+-}\circ \psi^{-+}.
\end{equation}
Equations (\ref{incl}) imply that there exists $\psi^{+-}$ defined on $\ker(D_X^{+-})$ and $\psi^{-+}$
 defined on $\ker(D_X^{+-})$ such that $\phi^{++}=D_Y^{-+}\circ \psi^{+-}$ and $\phi^{--}=-D_Y^{+-}\circ \psi^{-+}$. 
 Therefore, equations (\ref{WWW}) are satisfied on  $\ker(D_X^{+-})$ and  $\ker(D_X^{+-})$ 
 respectively. 
 
 \bigskip 
 \begin{rmk} Notice that any choice of  $\psi^{+-}$  on $\ker(D_X^{+-})$ and $\psi^{-+}$
  on $\ker(D_X^{+-})$ satisfying  $\phi^{++}=D_Y^{-+}\circ \psi^{+-}$ and $\phi^{--}=-D_Y^{+-}\circ \psi^{-+}$
  can be modified by adding to $\psi^{+-}$ any morphism from $\ker (D_X^{+-})$ to $\ker(D_Y^{-+})$ and to 
  $\psi^{-+}$ any morphism from $\ker (D_X^{-+})$ to $\ker(D_Y^{+-})$ . 
 \end{rmk}
 \bigskip 
 
 The problem is to extend $\psi^{+-}$ and  $\psi^{-+}$ to $X\otimes S^+$ and $X\otimes S^-$
 respectively. Since $\ker (S)\subset \ker(TS)$, we see from the description of $\ker(TS)$ given above
 that $\phi$ is $0$ on the $\widetilde K$-isotypic components which are not in the kernel of $D^2$.
 Therefore, we may set $\psi^{+-}$ and  $\psi^{-+}$ to be $0$ on these components, and 
 (\ref{WWW}) will be satisfied on them. So it remains to define  $\psi^{+-}$  (resp.  $\psi^{-+}$)
 on some complement $W^+$  of $\ker( D_X^{+-})$ in $\ker (D_X^{-+}\circ D_X^{+-})$
(resp.  some complement $W^-$ of $\ker( D_X^{-+})$ in $\ker (D_X^{+-}\circ D_X^{-+})$).

Let $x\in W^+$. Then $D_X^{+-}(x)\in \ker D_X ^{-+}$, so $\psi^{-+}(D_X^{+-}(x))$ is defined, and 
\[  \phi^{--}(  D_X^{+-}(x))= -D_Y^{-+} ( \psi^{-+}(D_X^{+-}(x))) .\]
But since $\phi \in \ker S$, Equations (\ref{kerS}) hold, and thus 
\[  \phi^{--}(  D_X^{+-}(x))=D_Y^{+-}(\phi^{++}(x))= -D_Y^{-+} ( \psi^{-+}(D_X^{+-}(x))) ,\]
and $\phi^{++}(x)+\psi^{-+}(D_X^{+-}(x))\in \ker(D_Y^{-+})$. Notice that $D_X^{+-}$ induces an isomorphism
from $W^-$ to $\ker (D_X^{-+})\cap \im (D_X^{+-})$. But by the remark above, we may assume that 
in fact $\phi^{++}+\psi^{-+}\circ D_X^{+-}$ is $0$ on $W^+$. Setting $\psi^{+-}$ to be $0$ on $W^+$, we see that 
  the first equation of  (\ref{WWW}) is satisfied on $W^+$. Similarly, we extend $\psi^{-+}$.
This finishes the proof of the injectivity of the first morphism in Lemma \ref{chiant}.
\medskip 

We now prove the surjectivity. Let : 
\[ \alpha : \, \frac{\ker (D_X^{+-}) }{ \ker (D_X^{+-})\cap  \im(D_X^{-+})  }
\longrightarrow \frac{\ker (D_Y^{+-}) }{\ker (D_Y^{+-})\cap  \im(D_Y^{-+})} , \]
and 
\[\beta : \, \frac{\ker (D_X^{-+}) }{ \ker (D_X^{-+})\cap  \im(D_X^{+-})}  
\longrightarrow \frac{\ker (D_Y^{-+}) }{ \ker (D_Y^{-+})\cap  \im(D_Y^{+-})} . \]
We would like to find $\phi=\phi^{++}+\phi^{--}$ such that $(\bar \phi^{++}, \bar\phi^{--})=(\alpha,\beta)$. 
First,  lift $\alpha$ to a morphism $\tilde \alpha : \, \ker(D_X^{+-})\rightarrow \ker(D_Y^{-+})$ 
which is identically $0$ on $\ker (D_X^{+-})\cap  \im(D_X^{-+})  $ and set $\phi^{++}=\tilde \alpha$ on 
$\ker(D_X^{+-})$.
Since $\phi$ should be in $\ker S$, we set $\phi$ to be $0$ on all 
$\widetilde K$-isotypic components not in the kernel of $D^2$. Let $W^+$, $W^-$  be as above.
Then, it remains to define $\phi^{++}$ on $W^+$ and we do this by setting  $\phi^{++}=0$ on $W^+$.
Similarly, extend $\phi^{--}$. It is then  immediately clear that Equations (\ref{kerS})
are satisfied, so we have constructed $\phi$ as we wished. This finishes the proof of Lemma \ref{chiant}.
\qed
\bigskip

Let us introduce now the index $\ind(S,T)$ of the Fredholm pair $(S,T)$. 
The material about Fredholm pairs and their indices is exposed in the appendix, and by definition
\[ \ind(S,T)=     \dim \left( \frac{\ker S}{ \ker S\cap \im T}\right) 
-\dim\left( \frac{\ker T}{ \ker T\cap \im S }\right)  . \]

\bigskip

\begin{cor}\[\ind(S,T)=\left[ \Idir(X) , \Idir(Y)\right]_{\widetilde K} \]
\end{cor}
\proof 
\begin{align*} \ind(S,T)&=     \dim \left( \frac{\ker S}{ \ker S\cap \im T}\right) -\dim\left( \frac{\ker T}{ \ker T\cap \im S }\right)    \\
&=\dim \left( \Hom_{\widetilde K}  \left( \Hdir(X)^+,\Hdir(Y)^+\right)\right)  +
\dim \left( \Hom_{\widetilde K}  \left(  \Hdir(X)^-  ,\Hdir(Y)^-  \right)\right)\\
&- 
\dim \left(  \Hom_{\widetilde K}  \left( \Hdir(X)^+,\Hdir(Y)^-\right) \right)  -
\dim \left( \Hom_{\widetilde K}  \left(  \Hdir(X)^-  ,\Hdir(Y)^+  \right)\right).\\
&=\left[ \Idir(X) , \Idir(Y)\right]_{\widetilde K} \end{align*}

Conjecture  \ref{Main1} is then equivalent to the fact that $\ind(S,T)$ is equal
to the Euler-Poincaré characteristic of the complex computing the $\Ext$ groups of $X$ and $Y$, namely
$\mathbf C=\sum_{i\in \bbN} C^i$, where $C^i=\Hom_K( \bigwedge^i\frp\otimes X,Y)$.
By example \ref{Fredholmcomplexe}, this Euler-Poincaré characteristic is 
also the index of the Fredholm pair given by the differentials between the even and odd part of the complex.
More precisely, with 
\[d^+: \mathbf{C}^{\bar 0} =\bigoplus_{i\in \bbN} C^{2i}\longrightarrow \mathbf{C}^{\bar 1}=
 \bigoplus_{i\in \bbN} C^{2i+1},\]
 and $d^-:  \mathbf{C}^{\bar 1} \rightarrow \mathbf{C}^{\bar 0}$, we have 
 $\EP(X,Y)=\ind(d^+, d^-)$.
So we would like to show that $\ind(S,T)=\ind(d^+, d^-)$.  To facilitate the comparison, first notice that 
\[ \mathbf C=\bigoplus_{i\in \bbN} C^i = \bigoplus_{i\in \bbN} \Hom_K({\textstyle  \bigwedge}^i\frp\otimes X,Y)
\simeq \Hom_K( \bigwedge\frp\otimes X,Y)\simeq \Hom_{\widetilde K}(X\otimes S, Y\otimes S)=\scrC \]
Then transport $(S,T)$ to $(\scrD^+, \scrD^-)$ via this isomorphism. Thus 
\[\scrD^+:\, \mathbf{C}^{\bar 0}\longrightarrow \mathbf{C}^{\bar 1}, 
\quad \scrD^-:\, \mathbf{C}^{\bar 1}\longrightarrow \mathbf{C}^{\bar 0}, \]
and set $\scrD=\scrD^+\oplus \scrD^-$, an operator on $\mathbf C$.
The conjecture is now that 
\begin{equation}\label{Conj2}\ind (d^+, d^-)= \ind(\scrD^+, \scrD^- ).
\end{equation}
\bigskip 

We give now explicit formulas for $d$ and $\scrD$.
The differential $d$ on $\mathbf C$ is the sum of the $d^i: C^i \rightarrow C^{i+1}$
\begin{align}
&d^i\phi(X_0\wedge\ldots \wedge X_i\otimes x)=\\
\nonumber \sum_{j=0}^i (-1)^j \big( X_j\cdot & \phi
(X_0\wedge\ldots\wedge \hat X_j \wedge\ldots \wedge  X_i\otimes x) 
- \phi(X_0\wedge\ldots\wedge \hat X_j \wedge\ldots \wedge  X_i\otimes X_j\cdot x) \big)
\end{align}

Recall that $\frp=U\oplus U^*$, where $U$ and $U^*$ are maximal isotropic subspaces
for the bilinear form $B$, and $U^*$, as the notation suggests, is identified with the dual of 
$U$ via $B$. We have realized the Clifford module $S$ as $\bigwedge U$, and thus  $S^*$ as 
$\bigwedge U^*$. Let us denote respectively $\gamma$ and $\gamma^*$ the Clifford actions of $\Cl(\frp)$ on $S$ and $S^*$.
Explicitely, for $u\in U$, $u^*\in U^*$, $\lambda_1\wedge \ldots \wedge \lambda_r\in \bigwedge^r  U \subset S$, and 
$\mu_1\wedge \ldots \wedge \mu_s\in \bigwedge^s  U^*\subset S^*$,
\begin{align}
&\gamma(u)(\lambda_1\wedge \ldots \wedge \lambda_r)=u\wedge\lambda_1\wedge \ldots \wedge \lambda_r\\
&\gamma(u^*)(\lambda_1\wedge \ldots \wedge \lambda_r)=2\sum_{j=1}^r (-1)^j B(u^*,\lambda_j)\; 
\lambda_1\wedge \ldots    \wedge \hat \lambda_j \wedge \ldots   \wedge \lambda_r\\
&\gamma^*(u^*)(\mu_1\wedge \ldots \wedge \mu_s)=u^*\wedge \mu_1\wedge \ldots \wedge \mu_s \\
 &\gamma^*(u)(\lambda_1\wedge \ldots \wedge \lambda_r)=2\sum_{j=1}^r (-1)^j B(u,\mu_j)\; 
\mu_1\wedge \ldots    \wedge \hat \mu_j \wedge \ldots   \wedge \mu_s
\end{align}

Fix a basis $u_1,\ldots ,u_m$ of $U$, with dual basis $u_1^*, \ldots ,u_m^*$. 
 Decompose the Dirac operator as 
 \[D= A+B\in U(\frg)\otimes \Cl(\frp), \quad A=\sum_{i=1}^m  u_i \otimes u_i^*, \quad B=\sum_{i=1}^m u_i^* \otimes u_i  \]
Then $A$ acts on $\mathbf{C}= \Hom_K(\bigwedge \frp\otimes X,Y)\simeq \Hom_K(S\otimes S^*\otimes X,Y)$
as follows :  for  $\lambda_1\wedge \ldots \wedge \lambda_r\in \bigwedge^r   U\subset S$, 
$\mu_1\wedge \ldots \wedge \mu_s\in \bigwedge^s U^*\subset S^*$, $x\in X$ and $\phi \in \Hom_K(S\otimes S^*\otimes X,Y)$, 
\begin{align} \label{A}
&(A\cdot \phi ) (  \lambda_1\wedge \ldots \wedge \lambda_r\wedge \mu_1\wedge \ldots \wedge \mu_s \otimes x)\\
\nonumber 
=\sum_{i=1}^m  
& \big[ u_i\cdot \phi( \big(\gamma(u_i^*)(\lambda_1\wedge \ldots \wedge \lambda_r)\big)
\wedge (\mu_1\wedge \ldots \wedge \mu_s)  \otimes  x)\\
\nonumber &+ (-1)^r\;  u_i\cdot  \phi( (\lambda_1\wedge \ldots \wedge \lambda_r)
\wedge \left( \gamma^*(u_i^*)( \mu_1\wedge \ldots \wedge \mu_s)\right)  \otimes  x) \\
\nonumber &- \phi( \left( \gamma(u_i^*)(\lambda_1\wedge  \ldots \wedge \lambda_r)\right) 
\wedge (\mu_1\wedge \ldots \wedge \mu_s) \otimes u_i \cdot x) \\
\nonumber &- (-1)^r \;  \phi( (\lambda_1\wedge  \ldots \wedge \lambda_r)
 \wedge\left( \gamma^*(u_i^*)(\mu_1\wedge \ldots \wedge \mu_s)\right)  \otimes u_i\cdot  x)\big]\\
\nonumber  = 2 \sum_{i=1}^m   & \sum_{j=1}^r (-1)^j \;     B(u_i^*,\lambda_j)\; 
u_i\cdot \phi( \lambda_1\wedge \ldots \wedge \hat \lambda_j \wedge \ldots  \wedge \lambda_r
\wedge (\mu_1\wedge \ldots \wedge \mu_s)  \otimes  x)\\
\nonumber+\sum_{i=1}^m & (-1)^r \;  u_i\cdot  \phi( (\lambda_1\wedge \ldots \wedge \lambda_r)
\wedge (u_i^*\wedge  \mu_1\wedge \ldots \wedge \mu_s)  \otimes  x) \\
\nonumber-  2 \sum_{i=1}^m   & \sum_{j=1}^r (-1)^j   \;  B(u_i^*,\lambda_j)\; 
\phi( \lambda_1\wedge \ldots \wedge \hat \lambda_j \wedge \ldots  \wedge \lambda_r
\wedge (\mu_1\wedge \ldots \wedge \mu_s)  \otimes u_i\cdot  x)\\
\nonumber-\sum_{i=1}^m & (-1)^r  \;  \phi( (\lambda_1\wedge \ldots \wedge \lambda_r)
\wedge (u_i^*\wedge  \mu_1\wedge \ldots \wedge \mu_s)  \otimes u_i\cdot  x) \\
\nonumber  = 2    \sum_{j=1}^r & (-1)^j\;      \, 
\lambda_j \cdot \phi( \lambda_1\wedge \ldots \wedge \hat \lambda_j \wedge \ldots  \wedge \lambda_r
\wedge (\mu_1\wedge \ldots \wedge \mu_s)  \otimes  x)\\
\nonumber+\sum_{i=1}^m & (-1)^r \;   u_i\cdot  \phi( (\lambda_1\wedge \ldots \wedge \lambda_r)
\wedge (u_i^*\wedge  \mu_1\wedge \ldots \wedge \mu_s)  \otimes  x) \\
\nonumber-  2     \sum_{j=1}^r & (-1)^j \;      
\phi( \lambda_1\wedge \ldots \wedge \hat \lambda_j \wedge \ldots  \wedge \lambda_r
\wedge (\mu_1\wedge \ldots \wedge \mu_s)  \otimes \lambda_j \cdot  x)\\
\nonumber-\sum_{i=1}^m & (-1)^r \;   \phi( (\lambda_1\wedge \ldots \wedge \lambda_r)
\wedge (u_i^*\wedge  \mu_1\wedge \ldots \wedge \mu_s)  \otimes u_i\cdot  x)  .
\end{align}
 
In this computation, notice that when we identify $\bigwedge \frp$ with $S\otimes S^*$, the 
latter is endowed with the super tensor product, and accordingly for the  Clifford
action. This explain the $(-1)^r$ factors in the formulas above. Similarly, the action of $B$ is given by 

\begin{align}\label{B} 
&(B\cdot \phi ) (  \lambda_1\wedge \ldots \wedge \lambda_r\wedge \mu_1\wedge \ldots \wedge \mu_s \otimes x)\\
\nonumber 
=\sum_{i=1}^m  
& \big[ u_i^*\cdot \phi( \big(\gamma(u_i)(\lambda_1\wedge \ldots \wedge \lambda_r)\big)
\wedge (\mu_1\wedge \ldots \wedge \mu_s)  \otimes  x)\\
\nonumber &+ (-1)^r\;  u_i^*\cdot  \phi( (\lambda_1\wedge \ldots \wedge \lambda_r)
\wedge \left( \gamma^*(u_i)( \mu_1\wedge \ldots \wedge \mu_s)\right)  \otimes  x) \\
\nonumber &- \phi( \left( \gamma(u_i)(\lambda_1\wedge  \ldots \wedge \lambda_r)\right) 
\wedge (\mu_1\wedge \ldots \wedge \mu_s) \otimes u_i^* \cdot x) \\
\nonumber &- (-1)^r \;  \phi( (\lambda_1\wedge  \ldots \wedge \lambda_r)
 \wedge\left( \gamma(u_i)(\mu_1\wedge \ldots \wedge \mu_s)\right)  \otimes u_i^*\cdot  x)\big]\\
\nonumber  =  \sum_{i=1}^m   & 
u_i^* \cdot \phi( (u_i\wedge \lambda_1\wedge \ldots  \wedge \lambda_r)
\wedge (\mu_1\wedge \ldots \wedge \mu_s)  \otimes  x)\\
\nonumber  + 2 (-1)^r  \sum_{i=1}^m   & \sum_{j=1}^s (-1)^j \;     B(u_i,\mu_j)\; 
u_i^* \cdot \phi( (\lambda_1\wedge \ldots  \wedge \lambda_r)
\wedge (\mu_1\wedge \ldots  \wedge \hat \mu_j \wedge \ldots    \wedge \mu_s)  \otimes  x)\\
\nonumber -\sum_{i=1}^m   & 
 \phi( (u_i\wedge \lambda_1\wedge \ldots  \wedge \lambda_r)
\wedge (\mu_1\wedge \ldots \wedge \mu_s)  \otimes  u_i^* \cdot  x)\\
\nonumber  - 2 (-1)^r  \sum_{i=1}^m   & \sum_{j=1}^s (-1)^j \;     B(u_i,\mu_j)\; 
 \phi( (\lambda_1\wedge \ldots  \wedge \lambda_r)
\wedge (\mu_1\wedge \ldots  \wedge \hat \mu_j \wedge \ldots    \wedge \mu_s)  \otimes  u_i^* \cdot x)\\
\nonumber  =  \sum_{i=1}^m   & 
u_i^* \cdot \phi( (u_i\wedge \lambda_1\wedge \ldots  \wedge \lambda_r)
\wedge (\mu_1\wedge \ldots \wedge \mu_s)  \otimes  x)\\
\nonumber  + 2 (-1)^r    & \sum_{j=1}^s (-1)^j \;     
\mu_j^* \cdot \phi( (\lambda_1\wedge \ldots  \wedge \lambda_r)
\wedge (\mu_1\wedge \ldots  \wedge \hat \mu_j \wedge \ldots    \wedge \mu_s)  \otimes  x)\\
\nonumber -\sum_{i=1}^m   & 
 \phi( (u_i\wedge \lambda_1\wedge \ldots  \wedge \lambda_r)
\wedge (\mu_1\wedge \ldots \wedge \mu_s)  \otimes  u_i^* \cdot  x)\\
\nonumber  - 2 (-1)^r     & \sum_{j=1}^s (-1)^j \;     
 \phi( (\lambda_1\wedge \ldots  \wedge \lambda_r)
\wedge (\mu_1\wedge \ldots  \wedge \hat \mu_j \wedge \ldots    \wedge \mu_s)  \otimes  \mu_j^* \cdot x).
\end{align}

The action of $D$ is the sum of the actions of $A$ and $B$.
Consider the four  operators defined by the last four lines of (\ref{A}) and denote them respectively by 
\[  2d_1=\scrE_1, \quad \delta_1=\scrD_1, \quad 2d_2=\scrD_2, \quad \delta_2=\scrE_2\]
and similarly 
for the four operators defined by the last four lines of (\ref{B}) 
\[ \delta_3=\scrE_3, \quad 2d_3=\scrD_3, \quad \delta_4=\scrD_4, \quad 2d_4=\scrE_4. \]

 Then we have 
 \[d=d_1+d_2+d_3+d_4, \qquad \scrD=\scrD_1+\scrD_2+\scrD_3+\scrD_4.\]
 (actually, the operators $S$ and $T$  of $\scrC$ at
 the begining of this section where  defined to make this true).

We get easily now  that $\partial=\partial_1+\partial_2+\partial_3+\partial_4$ satisfies $\partial^2=0$.
It is likely that the cohomology of $\partial$ in degree $i$ comptutes $\Ext^i(Y\, \check{},X\, \check{}\, )$
$Y\, \check{}$ and $X\, \check{}$ being the contragredient Harish-Chandra modules of $Y$ and $X$ respectively,
but we haven't checked this. The operator  $\scrE=\scrE_1+\scrE_2+\scrE_3+\scrE_4$
is very similar to $\scrD$. It is obtained as $\scrD$ from the transport of $(T,S)$, but using a different isomorphism 
\[ \mathbf C=\Hom_K( \bigwedge\frp\otimes X,Y)\simeq \Hom_{\widetilde K}(X\otimes S, Y\otimes S)=\scrC \]
obrtained by exchanging the role of the two copies of $S$ after identifying it with its dual.
Thus $d$ and $\scrD$ are related by 
\begin{equation}\label{BigOp}2d+\delta=\scrD+\scrE.\end{equation}
A possible way to prove (\ref{Conj2}) would be to show that both indices are equal to the index of the 
operator (\ref{BigOp}).

\section{Appendix : Fredholm pairs}\label{Fred}

In this section, we adapt from \cite{Amb} the definition of the index of a Fredholm pair, as well as some properties
of this invariant which enable one to calculate it in pratical applications. We do this in a purely algebraic setting, while
the theory is developed for Banach spaces in \cite{Amb}.

\begin{defi}\label{FredholmPair} Let $X$, $Y$ be complex vector spaces and let $S \in  \caL(X; Y )$, $T \in \caL(Y;X)$.
Then $(S, T )$ is called a Fredholm pair if the following dimensions are finite :
\[ a := \dim  \ker(S)/ \ker(S) \cap  \im(T );    \quad b := \dim \ker (T) /  \ker(T ) \cap \im(S).\]
In this case, the number
\[  \ind(S, T) := a -b \]
is called the index of $(S, T )$.
\end{defi}

\bigskip

\begin{example}\label{FredhomOp}
Take $T=0$. Then $(S,0)$ is a Fredholm pair if and only if $S$ is a Fredholm operator
and 
\[\ind(S,0)=\ind(S)=\dim\ker (S)-\dim\coker (S)\]
\end{example}

\begin{example}\label{Fredholmcomplexe}
Consider a differential complex
\[ \cdots\longrightarrow  C^{i-1}\stackrel{d^{i-1}}{\longrightarrow } C^{i}\stackrel{d^{i}}{\longrightarrow }C^{i+1}
\longrightarrow \cdots \]
Suppose the cohomology groups $H^i:=\ker d^i/ \im d^{i-1}$ of this complex are finite dimensional, and non zero only 
for a finite number of them.
Put $X=\bigoplus_{i \in \bbZ}C^{2i}$, $Y=\bigoplus_{i \in \bbZ}C^{2i+1}$, $S=\oplus_{i \in 2\bbZ} d^{2i}$, 
$T=\oplus_{i \in 2\bbZ} d^{2i+1}$. Then $(S,T)$ is a Fredholm pair and its index is equal to the Euler-Poincaré
characteristic of the complex : 
\[\ind(S,T)=\sum_{i\in \bbZ} (-1)^i  \dim H^i .\]
\end{example}


\bigskip 

With the notation of the definition, notice that $T$ induces an isomorphism 
\[\widetilde T: \, \frac{\im (S)}{\im(S)\cap \ker(T)} \longrightarrow \im(TS)\] 
and that $S$  induces an isomorphism 
\[ \widetilde S:\,  \frac{\im (T)}{\im(T)\cap \ker(S)}\longrightarrow \im(ST).\] 
Since $S(\im(TS))=\im(STS)\subset \im (ST)$, we see also that
$S$ induces a morphism
\[ \bar S: \, X/ \im(TS)\longrightarrow Y /\im(ST).\]
Similarly $T$ induces a morphism
\[\bar T: \, Y/\im(ST)\longrightarrow X /\im(TS).\]

Set $\bar X= X/\im(TS)$, $\bar Y=Y /\im(ST)$. 
\begin{lemma}\label{tools}
Suppose that $(S, T)$ is a Fredholm pair. Then $(\bar S ,  \bar T)$ is also a Fredholm pair and 
$\ind(S,T)=\ind(\bar S,\bar T)$. Suppose furthermore than $\bar X$ and $\bar Y$ are finite dimensional, 
then $\ind(\bar S, \bar T)=\dim \bar X-\dim(\bar Y)$. 
\end{lemma}

We leave the easy proof to the reader. Notice that $\bar T \bar S=0$ and that $\bar S\bar T=0$, so that 
$\im(\bar S)\subset \ker (\bar T)$, $\im(\bar T)\subset \ker (\bar S)$.

\begin{lemma}\label{exactsequ}
Suppose $(S_j,T_j)$ are pairs of linear maps between spaces $X_j$ and $Y_j$  with $S_jT_j=0$, 
$T_jS_j=0$, $j=1,2,3$.
Suppose that the following diagram is commutative, with vertical exact sequences : 
\[ \xymatrix{
0 \ar[d] & 0  \ar[d]& 0 \ar[d] \\
 X_1 \ar[r]^{S_1} \ar[d]^{\alpha}& Y_1 \ar[r]^{T_1} \ar[d]^{\gamma} & X_1 \ar[d]^{\alpha}\\
  X_2 \ar[r]^{S_2}\ar[d]^{\beta}& Y_2 \ar[r]^{T_2}\ar[d]^{\delta}& X_2\ar[d]^{\beta}\\ 
  X_3\ar[r]^{S_3}\ar[d] & Y_3\ar[r]^{T_3} \ar[d]& X_3\ar[d]\\
  0 & 0 & 0
}.\]

Then if  two of the pairs $(S_j; T_j)$ are Fredholm, the third is, and 
\[ \ind(S_1,T_1)-\ind(S_2,T_2)+\ind(S_3,T_3)=0.\]
\end{lemma}

This is \cite{Amb}, Lemma 2.2.

\bigskip
Consider now the following situation. Suppose $V=V^{\bar 0}\oplus V^{\bar 1}$ is a super vector space.
Suppose that $d$, $\partial$ are odd operators on $V$,  satisfying $d^2=\partial^2=0$.
Set  $d^{+}, \partial^{+} : V^{\bar 0 } \rightarrow V^{\bar 1} $    and  
$d^{-}, \partial^{-} : V^{\bar 1 } \rightarrow V^{\bar 0} $ for the restriction of $d$, $\partial$
to the even and odd part respectively. Suppose that the cohomology groups
for $d$, $H_d^+=\ker d^+/\im(d^-)$ and $H_d^-=\ker d^-/\im(d^+)$ are finite dimensional. 
Then according to example \ref{Fredholmcomplexe}, $(d^+;d^-)$
is a Fredholm pair and $\ind(d^+, d^-)=\dim H_d^+-H_d^-$.

Set $\scrF=d+\partial$. This is an odd operator on $V$ and similarly to $d$ and $\partial$, we denote
by $\scrF^{+}$ and $\scrF^{-}$ its restriction to the even and odd part of $V$ respectively.
The following result could be useful in proving (\ref{Conj2}) via (\ref{BigOp}).

\begin{prop}\label{mainFred}  Suppose $\scrF$ has the following properties : 
$\ker (\scrF^2)$ is finite dimensional and $\ker (\scrF^2)\oplus \im(\scrF^2)=V$.  

Then  $(\scrF^{+}, \scrF^{-})$ is a Fredholm pair and $\ind(\scrF^{+}, \scrF^{-})=\ind(d^{+}, d^{-})$.
\end{prop}

\proof The first assertion is obvious. From Lemma \ref{tools}, we have furthermore than 
$\ind(\scrF^{+}, \scrF^{-})  = \dim( (\ker \scrF^2)^{\bar 0} )- \dim ((\ker \scrF^2)^{\bar 1} )$, where we have set 
\[ (\ker \scrF^2)^{\bar 0}=\ker (\scrF^2)\cap V^{\bar 0}=\ker(\scrF^-\scrF^+), \quad 
(\ker \scrF^2)^{\bar 1}=\ker (\scrF^2)\cap V^{\bar 1}=\ker(\scrF^+\scrF^-). \]
 
Let us now compute the index of the pair $(d^+,d^-)$. First note that $\im \scrF^2$ is stable under $d$.
Indeed, if $x=  \scrF^2(y)= (d\partial +\partial d)(y)\in \im (\scrF^2)$, then 
$d(x)= d^2\partial (y)+d\partial d(y)=d\partial d(y)=(d\partial+\partial d )(d(y))=\scrF^2(dy)\in \im (\scrF^2)$. 
Thus $d$ induces operators, that we will still denote by $d$, on the spaces $\im(\scrF^2)$ and $V/\im(\scrF^2)$.
These spaces inherits from $V$ the structure of super vector spaces, and thus operators $d^+$, $d^-$
are well defined on them with the obvious meaning. Let us apply Lemma \ref{exactsequ} to the exact sequence 
\[ 0\longrightarrow \im (\scrF^2)\longrightarrow V\longrightarrow V/\im(\scrF^2)\longrightarrow 0.\]
Set $W=\im(\scrF^2)$, and let us compute the index of $(d^+,d^-)$ acting on $W=W^{\bar 0}\oplus  W^{\bar 1}$.
We claim  that : 
\begin{equation}W=\ker \partial+\im d
\end{equation}
Indeed, let $x=  \scrF^2(y)= (d\partial +\partial d)(y)\in \im (\scrF^2)$. Then $\partial(x)=\partial d\partial (y)$, 
so $x-d \partial (y)\in \ker \partial$.  This proves the claim, since $d\partial (y)\in \im d$.

Now, suppose that $x\in W$ is in $\ker d$.  Write $x=z+d(y)$, with $z\in \ker \partial$,  as we have just shown possible. 
Then $0=d(x)=d(z)$ so $z\in \ker d$. But $\scrF^2=d\partial+\partial d$ implies that
$\ker d\cap \ker \partial\subset \ker (\scrF^2)$. But by assumption,  
$\im(\scrF^2)\cap \ker(\scrF^2)=\{0\}$, and thus $z=0$. This shows that $x=d(y)\in \im d$.
Thus, the index  of $(d^+,d^-)$ acting on $W$ is $0$.  The Lemma  \ref{exactsequ} then implies
that  the index  of $(d^+,d^-)$ acting on $W$ is equal to the index of $(d^+,d^-)$ acting on $V/W$.
But by the semisimplicity of $\scrF^2$, $V/W$ is isomorphic to $\ker (\scrF^2)$ as super vectors spaces, 
and therefore finite dimensional. Another application of Lemma \ref{tools} then tells us that 
the index of  $(d^+,d^-)$ acting on $V/W$ equals $\dim( (\ker \scrF^2)^{\bar 0} )- \dim ((\ker \scrF^2)^{\bar 1} )$.
\qed

\bibliographystyle{abbrv}
\bibliography{EP}

\def\cftil#1{\ifmmode\setbox7\hbox{$\accent"5E#1$}\else
  \setbox7\hbox{\accent"5E#1}\penalty 10000\relax\fi\raise 1\ht7
  \hbox{\lower1.15ex\hbox to 1\wd7{\hss\accent"7E\hss}}\penalty 10000
  \hskip-1\wd7\penalty 10000\box7}
\begin{thebibliography}{10}

\bibitem{Amb}
C.-G. Ambrozie.
\newblock On {F}redholm index in {B}anach spaces.
\newblock {\em Integral Equations Operator Theory}, 25(1):1--34, 1996.

\bibitem{Ar}
J.~Arthur.
\newblock On elliptic tempered characters.
\newblock {\em Acta Math.}, 171(1):73--138, 1993.

\bibitem{BlBr}
P.~Blanc and J.-L. Brylinski.
\newblock Cyclic homology and the {S}elberg principle.
\newblock {\em J. Funct. Anal.}, 109(2):289--330, 1992.

\bibitem{Bouaziz}
A.~Bouaziz.
\newblock Int\'egrales orbitales sur les groupes de {L}ie r\'eductifs.
\newblock {\em Ann. Sci. \'Ecole Norm. Sup. (4)}, 27(5):573--609, 1994.

\bibitem{Bouaziz2}
A.~Bouaziz.
\newblock Formule d'inversion des int\'egrales orbitales sur les groupes de
  {L}ie r\'eductifs.
\newblock {\em J. Funct. Anal.}, 134(1):100--182, 1995.

\bibitem{COT}
D.~Ciubotaru, E.~Opdam, and P.~Trapa.
\newblock Algebraic and analytic dirac induction for graded affine hecke
  algebras.
\newblock {\em Arxiv arXiv:1201.2130}.

\bibitem{CT}
D.~M. Ciubotaru and P.~E. Trapa.
\newblock Characters of {S}pringer representations on elliptic conjugacy
  classes.
\newblock {\em Duke Math. J.}, 162(2):201--223, 2013.

\bibitem{Dat}
J.-F. Dat.
\newblock On the {$K_0$} of a {$p$}-adic group.
\newblock {\em Invent. Math.}, 140(1):171--226, 2000.

\bibitem{Dat2}
J.-F. Dat.
\newblock Une preuve courte du principe de {S}elberg pour un groupe
  {$p$}-adique.
\newblock {\em Proc. Amer. Math. Soc.}, 129(4):1213--1217, 2001.

\bibitem{HC}
Harish-Chandra.
\newblock Supertempered distributions on real reductive groups.
\newblock In {\em Studies in applied mathematics}, volume~8 of {\em Adv. Math.
  Suppl. Stud.}, pages 139--153. Academic Press, New York, 1983.

\bibitem{HKP}
J.-S. Huang, Y.-F. Kang, and P.~Pand{\v{z}}i{\'c}.
\newblock Dirac cohomology of some {H}arish-{C}handra modules.
\newblock {\em Transform. Groups}, 14(1):163--173, 2009.

\bibitem{HP}
J.-S. Huang and P.~Pand{\v{z}}i{\'c}.
\newblock Dirac cohomology, unitary representations and a proof of a conjecture
  of {V}ogan.
\newblock {\em J. Amer. Math. Soc.}, 15(1):185--202 (electronic), 2002.

\bibitem{HPbook}
J.-S. Huang and P.~Pand{\v{z}}i{\'c}.
\newblock {\em Dirac operators in representation theory}.
\newblock Mathematics: Theory \& Applications. Birkh\"auser Boston, Inc.,
  Boston, MA, 2006.

\bibitem{Kaz}
D.~Kazhdan.
\newblock Cuspidal geometry of {$p$}-adic groups.
\newblock {\em J. Analyse Math.}, 47:1--36, 1986.

\bibitem{KV}
A.~W. Knapp and D.~A. Vogan, Jr.
\newblock {\em Cohomological induction and unitary representations}, volume~45
  of {\em Princeton Mathematical Series}.
\newblock Princeton University Press, Princeton, NJ, 1995.

\bibitem{Kos1}
B.~Kostant.
\newblock Clifford algebra analogue of the {H}opf-{K}oszul-{S}amelson theorem,
  the {$\rho$}-decomposition {$C(\mathfrak g)={\rm End}\, V_\rho \otimes
  C(P)$}, and the {$\mathfrak g$}-module structure of {$\bigwedge \mathfrak
  g$}.
\newblock {\em Adv. Math.}, 125(2):275--350, 1997.

\bibitem{Kot}
R.~E. Kottwitz.
\newblock Tamagawa numbers.
\newblock {\em Ann. of Math. (2)}, 127(3):629--646, 1988.

\bibitem{Lab}
J.-P. Labesse.
\newblock Pseudo-coefficients tr\`es cuspidaux et {$K$}-th\'eorie.
\newblock {\em Math. Ann.}, 291(4):607--616, 1991.

\bibitem{Mein}
E.~Meinrenken.
\newblock {\em Clifford algebras and {L}ie theory}, volume~58 of {\em
  Ergebnisse der Mathematik und ihrer Grenzgebiete. 3. Folge. A Series of
  Modern Surveys in Mathematics [Results in Mathematics and Related Areas. 3rd
  Series. A Series of Modern Surveys in Mathematics]}.
\newblock Springer, Heidelberg, 2013.

\bibitem{Par}
R.~Parthasarathy.
\newblock Dirac operator and the discrete series.
\newblock {\em Ann. of Math. (2)}, 96:1--30, 1972.

\bibitem{ScSt}
P.~Schneider and U.~Stuhler.
\newblock Representation theory and sheaves on the {B}ruhat-{T}its building.
\newblock {\em Inst. Hautes \'Etudes Sci. Publ. Math.}, (85):97--191, 1997.

\bibitem{Vig2}
M.-F. Vign{\'e}ras.
\newblock Caract\'erisation des int\'egrales orbitales sur un groupe r\'eductif
  {$p$}-adique.
\newblock {\em J. Fac. Sci. Univ. Tokyo Sect. IA Math.}, 28(3):945--961 (1982),
  1981.

\bibitem{Vig}
M.-F. Vign{\'e}ras.
\newblock On formal dimensions for reductive {$p$}-adic groups.
\newblock In {\em Festschrift in honor of {I}. {I}. {P}iatetski-{S}hapiro on
  the occasion of his sixtieth birthday, {P}art {I} ({R}amat {A}viv, 1989)},
  volume~2 of {\em Israel Math. Conf. Proc.}, pages 225--266. Weizmann,
  Jerusalem, 1990.

\bibitem{Wal}
N.~R. Wallach.
\newblock {\em Real reductive groups. {I}}, volume 132 of {\em Pure and Applied
  Mathematics}.
\newblock Academic Press, Inc., Boston, MA, 1988.

\end{thebibliography}

\end{document}